\documentclass[a4paper, 10pt]{amsproc}
\usepackage{defs}
\usepackage{orcidlink}
\usepackage{calc}
\definecolor{metablue}{HTML}{0064E0}
\definecolor{metafg}{HTML}{1C2B33}
\definecolor{metabg}{HTML}{F1F4F7}
\definecolor{metabgdeep}{HTML}{D9EFFF}
\definecolor{metagreen}{HTML}{EAFFE8}
\definecolor{metagreen}{HTML}{FCFFEE}
\definecolor{metared}{HTML}{FFEAE8}
\usepackage[framemethod=tikz,xcolor=true]{mdframed}

\newmdenv[backgroundcolor=metabgdeep, roundcorner=10pt, skipabove=4pt, linewidth=0pt, innertopmargin=4pt]{myframe}

\newmdenv[backgroundcolor=metabg, roundcorner=10pt, skipabove=4pt, linewidth=0pt, innertopmargin=4pt]{myOCP}

\newmdenv[backgroundcolor=metared, roundcorner=10pt, skipabove=7pt, linewidth=0pt, innertopmargin=7pt]{myalgo}
\newmdenv[%
    leftmargin=0.5cm,
    backgroundcolor=yellow!10,%
    roundcorner=5pt,%
    tikzsetting={draw=red, line width=2.0pt}%
    ]{SpecialText}%

\title[RobHandsOff]{Robust maximum hands-off optimal control: existence, maximum principle, and \(\lpL[0]\)-\(\lpL[1]\) equivalence}

\author[S. Ganguly]{Siddhartha Ganguly\,\orcidlink{0000-0003-2046-2061}}
\author[K. Kashima]{Kenji Kashima \orcidlink{0000-0002-2963-2584}}
\thanks{%
	S. Ganguly and K. Kashima are with \faGroup\ The Applied Mathematics and Physics Department, Graduate School of Informatics, \faUniversity\ Kyoto University, \faMapMarker\  Kyoto, Japan. S.G. gratefully acknowledges Takuya Ikeda's valuable suggestions, as well as the encouraging comments from Debasish Chatterjee and Masaaki Nagahara.}

\thanks{%
	Contact Information: (SG) \faHome\ \url{https://sites.google.com/view/siddhartha-ganguly}, \faEnvelope\ \texttt{siddhartha1123@gmail.com},
    (KK): \faEnvelope\ \texttt{kk@i.kyoto-u.ac.jp}, \faHome\ \url{https://www.bode.amp.i.kyoto-u.ac.jp/en/kk/}. 
}

\begin{document}

\maketitle

\begin{abstract}
This work advances the maximum hands-off sparse control framework by developing a robust counterpart for constrained linear systems with parametric uncertainties. The resulting optimal control problem minimizes an \(\lpL[0]\) objective subject to an uncountable, compact family of constraints, and is therefore a nonconvex, nonsmooth robust optimization problem. To address this, we replace the \(\lpL[0]\) objective with its convex \(\lpL[1]\) surrogate and, using a nonsmooth variant of the robust Pontryagin maximum principle, show that the \(\lpL[0]\) and \(\lpL[1]\) formulations have identical sets of optimal solutions --- we call this \emph{the robust hands-off principle}. Building on this equivalence, we propose an algorithmic framework --- drawing on numerically viable techniques from the semi-infinite robust optimization literature --- to solve the resulting problems. An illustrative example is provided to demonstrate the effectiveness of the approach.
\end{abstract}

\section{Introduction}\label{sec:intro}

Sparsity, in broad strokes, is defined as the mathematical representation of \emph{simplicity}, capturing the principle of minimizing complexity in models or systems. Naturally, the study of sparsity and the associated numerical techniques to promote it have gained significant prominence across various scientific and technological fields, including signal processing, machine learning, and statistics, such as compressed sensing or sparse modeling \cite{ref:SparseBook:Unser,ref:SparseBook:Vidyasagar}, \cite{SparseBook:ML:StatLearning}, \cite{ref:SparseApp:Donoho}.

Sparsity-driven techniques have been extensively utilized in the control literature to minimize control activity across a wide range of applications. These include networked control systems \cite{ref:SparseApp:NetworkPackets}, feedback controller design \cite{ref:SparseApp:Feedgain, ref:PolKhlSch-14}, sensor placement and sensing strategies \cite{ref:SparseApp:SensingFeedDisgn}, aerospace systems \cite{ref:SparseApp:Space}, control of partial differential equations \cite{ref:SparseApp:PDEs:I, ref:SparseApp:PDE:II}, and infinite-dimensional systems \cite{ref:TI:MH:InfDim}. From a theoretical perspective, the system-theoretic properties of sparse optimal controllers constitute a highly active and rapidly developing area of research \cite{ref:agrachev2025optimal, ref:agrachev2025singularextremalsoptimalcontrol, ref:cavare2025l1}. By promoting sparse control inputs, these approaches optimize resource efficiency and reduce actuator engagement while maintaining system performance. From the perspective of constrained optimal control theory, the use of sparsity-promoting cost functions was explored in \cite{ref:MN:MaxHandsOff}. This approach underpins the well-established concept of the \highlight{maximum hands-off control} in the sparse optimal control literature. Within the context of deterministic linear systems, this principle is formalized by minimizing the \(\lpL[0]\)-(semi)norm of the control input, which effectively reduces the duration over which the control is nonzero \cite{ref:MN:MaxHandsOff}. By prioritizing sparsity, the maximum hands-off principle induces prolonged intervals of zero control input, thereby minimizing control effort by keeping actuators inactive for extended periods while ensuring that system constraints are satisfied. In practice, the nonconvex and nonsmooth \(\lpL[0]\) objective is often replaced by its convex surrogate, the \(\lpL[1]\) norm, a relaxation justified by theoretical results demonstrating their equivalence under mild assumptions on the problem data \cite{ref:MN:MaxHandsOff}. This substitution facilitates computational tractability while preserving the sparsity-promoting properties central to the maximum hands-off control framework, ultimately enhancing efficiency by reducing the active engagement of control actuators.


\textbf{Motivation, theme, and contributions:} Given the preceding premise, sparse control is typically not robust to uncertainties, and ensuring resilience to alterations in system models is a critical task. Real-world dynamical models contain uncertainties, including unmodeled dynamics and disturbances, which can degrade system performance if not included in the control design. The theoretical equivalence between the \(\lpL[0]\) and the \(\lpL[1]\) sparse optimal controls and their numerical synthesis is well-established \cite{ref:MN:MaxHandsOff} for deterministic uncertainty-free systems. On the otherhand, a rigorous theoretical foundation for robust sparse control has remained an open challenge, along with a numerically viable technique for the synthesis of such control. To this end, we ask:
\begin{myOCP}
   \begin{quote}
\emph{``Is it possible to establish an equivalence between \(\lpL[0]\)- and \(\lpL[1]\)-sparse controls for a class of parameter-dependent systems, and to develop a numerically viable method for synthesizing such control actions?''}
\end{quote} 
\end{myOCP}
This article provides an affirmative answer to this question; specifically, our contributions are as follows:
\begin{enumerate}[leftmargin=*]
\item We consider a class of constrained linear uncertain systems with the objective of minimizing the control activity, i.e., minimizing the \(\lpL[0]\) control cost, while satisfying an \emph{uncountable family} of constraints. The ensuing optimization problem is a nonconvex and nonsmooth robust program which is both challenging to analyze and solve in a numerically tractable fashion. 

\item Following the standard route \cite{ref:MN:MaxHandsOff}, we analyze its convex relaxation, the robust \(\lpL[1]\) problem. The main contribution of this work is a rigorous theoretical framework for this problem. By employing a nonsmooth version of the robust Pontryagin's Maximum Principle, we first derive the necessary conditions for optimality for the robust \(\lpL[1]\) problem. Using these conditions, we prove that under standard assumptions, the optimal control exhibits a \emph{bang-off-bang} structure. Leveraging this structural result, we then present the central theoretical finding of this paper: a proof that the set of optimal solutions for the original robust \(\lpL[0]\) problem is identical to that of its tractable \(\lpL[1]\) convex relaxation.

\item Finally, in the numerical front we follow the suggestions of \cite[\S I, p. 941]{ref:minmax:ness:vinter}: we parametrize the space of admissible controls employing piecewise constant functions to arrive at a finite dimensional semi-infinite optimization problem and then employ some recent algorithms to furnish a numerically tractable algorithmic architecture with guarantees of \emph{exact solutions}. We validate our theoretical claims and demonstrate the performance of the numerical architecture on a benchmark control example. 
\end{enumerate}

\textbf{Placement and potential applications of our results:} From a \emph{theoretical standpoint}, our results should be viewed as an addition to the literature on sparse optimal control and equivalence principles, filling a gap in the setting of robust parameter-dependent systems. While the route to equivalence is similar to standard works the technical and methodological engines are \emph{completely novel}. From an \emph{algorithmic viewpoint}, the proposed approach leverages existing robust optimization oracles; several such alternatives are available, including cutting-plane or cutting-surface methods \cite{ref:blankenship1976infinitely} and the scenario approach \cite{MC-SG:18} (see \cite{ref:Scenario:zhang2022sparse} for an application). We adopt the technique in \cite{ref:DasAraCheCha-22} as a representative method because of its ability to generate numerically viable \emph{exact solutions} in a computationally viable fashion, which is unique. Our results have potential applications in several areas, including attitude control of thrusters under uncertainty \cite{ref:kristiansen2023maximum}, networked control over communication channels subject to packet drops \cite{ref:SparseApp:NetworkPackets}, and energy-efficient control design under uncertainty \cite{ref:MN:23}. These directions will be investigated in future work.

The article is structured as follows: the central problem is formulated in \S\ref{sec:prob:form}. The primary theoretical contributions are presented in \S\ref{subsec:theo:dev:pmp}, while the algorithmic developments are detailed in \S\ref{subsec:algo:dev:sip}. Finally, \S\ref{sec:num_exp} provides a numerical example to illustrate our results and the algorithmic architecture.

\section{The central problem}\label{sec:prob:form}
\subsection*{Notation}
We employ standard notation in this article. We let \(\N \Let \aset{1,2,\ldots}\) denote the set of positive integers, and let \(\Z\) be the set of integers. Let \(d \in \N\); we will denote the standard Euclidean vector space by \(\Rbb^d\) which is assumed to be equipped with standard norm \(\Rbb^d \ni x\mapsto \norm{x}\). Let \(X\) be an arbitrary subset of \(\Rbb^d\); by \(\intr{X}\) we denote the interior of \(X\). Let \(\horizon>0\); for any bounded measurable function \(f:\lcrc{0}{\horizon}\lra \Rbb\) we employ the notation \(\norm{f(\cdot)}_{\infty} \Let \text{esssup}_{t\in\lcrc{0}{\horizon}}\abs{f(t)}\) for the uniform norm, and when \(f(\cdot)\) is continuous then the essential supremum `\text{esssup}' is replaced by the standard supremum `\(\sup\)' or the `max' operator. Let \(A\) be a subset of \(\Rbb\), then the indicator function of \(A\) is the function \(\indic{A}(\cdot)\) defined by \(\indic{A}(x) = 1\) if \(x\in A\) and \(\indic{A}(x) = 0\) otherwise. For a fixed \(\horizon > 0\) and fixed step \(h > 0\), the space of \(\Rbb\)-valued piecewise constant functions is given by
\begin{align}
   \pwcfunc_h(\lcrc{0}{\horizon};\Rbb) \Let \left\{\xi:\lcrc{0}{\horizon} \to \Rbb \;\middle\vert\;  
    \begin{array}{@{}l@{}}
		0=t_0 <t_1<\cdots<t_{K} \le \horizon,\; \xi_i\in\Rbb,\\
		\,t_i-t_{i-1} = h \text{ for }i=1,\ldots,K, T-t_{K} \le h,\\ 
		\xi(\cdot)  = \sum_{i=1}^{K} \xi_{i-1} \indic{\lcro{t_{i-1}}{t_i}}(\cdot) + \xi_{K}\indic{\lcrc{t_{K}}{\horizon}}(\cdot)
        \end{array}
        \right\}. \nn
\end{align}
Let \(d_1,d_2 \in \N\); for \(X\subset \Rbb^{d_1}\) and \(Y\subset \Rbb^{d_2}\) open subsets, the set of \(r\)-times continuously differentiable functions from \(X\) to \(Y\) is denoted by \(\mathcal{C}^r (X;Y)\) and for \(1 \le p < +\infty\) the space of Lebesgue measurable functions \(\lpL[p](X;Y)\) is defined by
\[
\lpL[p](X;Y) \Let \aset[\bigg]{u : X \lra Y \suchthat
u(\cdot) \text{ is measurable and} \int_X \norm{u(x)}^p\odif{x} < +\infty}.
\]
The set of measurable and essentially bounded functions \(\lpL[\infty](X;Y)\) is defined by
\[
\lpL[\infty](X;Y) \Let \aset[\bigg]{u : X \lra Y \suchthat
u(\cdot) \text{ is measurable and} \text{ esssup}_{x\in X}\norm{u(x)} < +\infty}.
\]
The space of absolutely continuous functions is defined by
\[
W^{1,1}(X;Y) \Let \aset[\bigg]{
u : X \lra Y \suchthat 
\begin{array}{l}
u(\cdot) \in \lpL[1](X;Y) \text{ and }\partial_j u_i(\cdot) \in \lpL[1](X) \\
\text{for all }i=1,\ldots,d_2, \text{ and }j=1,\ldots,d_1
\end{array}
}
\]
We assume that all these spaces are equipped with their standard norms \cite{ref:adams2003sobolev}. Finally, the normal cone to a given convex set \(C\) at a point \(x \in C\), denoted \(N_C(x)\), is defined by \cite[\S 1.4]{ref:clarke_ocpbook}
\begin{align}
   N_C(x) \Let  \aset[\big]{v \in \Rbb^d \suchthat \inprod{v}{y - x} \le 0 \text{ for all } y \in C}.\nn
\end{align}


We present the central problem of this manuscript. Let \(d,\nu \in \N\) and fix a time horizon \(\horizon>0\). Let us consider an uncertain linear time-invariant controlled dynamical system modeled by the ordinary differential equation
\begin{equation} \label{eq:sys}
\dot \st(t;\uparam) = A(\uparam) \st(t;\uparam) + b(\uparam) \cont(t)\, \text{for a.e.\, } t \in [\tinit, \tfin],\,\,\st(0;\uparam) \Let \param \text{ given}
\end{equation}
with the data: 
\begin{enumerate}[label={\highlight{\textup{(\eqref{eq:sys}-\alph*)}}}, leftmargin=*, widest=b, align=left]

\item \label{ocp:data:1} The uncertain parameter \(\uparam \in \pset\) where \(\pset \subset \Rbb^{\nu}\) is a nonempty, compact and convex set of \emph{uncertain parameters}. A commonly occurring example of such uncertainty set \(\pset\) is a polytope, and it is widely found in the robust control literature \cite[Chapter 4]{ref:BoyElGFerBal-94}, \cite[Chapter 8]{ref:DullPag-13}. 

\item \label{ocp:data:2} the control trajectory \(t\mapsto \cont(t) \in \Rbb\) is a bounded measurable function satisfying
\begin{equation}
    \label{eq:inst_control}
    \cont(t) \in \admcont \Let \lcrc{-1}{1} \subset \Rbb \quad\text{ for a.e. }t\in \lcrc{0}{\horizon}, \nn
\end{equation}
and consequently, the admissible set of control actions is defined by
\begin{align}
    \admcon \Let \aset[\big]{u(\cdot) \in \lpL[\infty](\lcrc{0}{\horizon};\Rbb) \suchthat \cont(t) \in \admcont \text{ for a.e. }t \in \lcrc{0}{\horizon}}. \nn
\end{align}
For every \(\uparam \in \pset\), the state trajectory \(\st(\cdot;\uparam) \in W^{1,1}(\lcrc{0}{\horizon};\Rbb^d)\) is an absolutely continuous function of time. We assume that the mappings \(\pset \ni \uparam \mapsto A(\uparam) \in \Rbb^{d \times d}\) and \(\pset \ni \uparam \mapsto  b(\uparam) \in \Rbb^{d}\) are continuously differentiable on an open neighbourhood of \(\pset\). We will call the pair \(\bigl(u(\cdot), \aset[]{\st(\cdot;\uparam) \suchthat \uparam \in \pset}\bigr)\) a \emph{process} satisfying the dynamics \eqref{eq:sys} with \(\st(0;\uparam) \Let \param\) given.

\item \label{ocp:data:3} The final state takes values in the compact and convex set \(C\) given by
\begin{align}
    C \Let  \aset[]{\xi \in \Rbb^d \suchthat \psi(\xi) \le 0}, \nn
\end{align}
i.e, \(\psi \bigl(\st(T;\uparam)\bigr) \le 0\) for all \(\uparam \in \pset\), where \(\psi: \Rbb^d \lra \Rbb\) is continuously differentiable and convex function.  
\end{enumerate}
\vspace{1mm}

The process \(\bigl(u(\cdot), \aset[]{\st(\cdot;\uparam) \suchthat \uparam \in \pset}\bigr)\) is \emph{feasible} if \(t \mapsto \cont(t) \in \admcont\) and \(\st(\horizon;\uparam) \in C\) for all \(\uparam \in \pset\). Define the support for the control trajectory by
\begin{align}
    \text{supp}(\cont) \Let \aset[\big]{t \in \lcrc{0}{\horizon} \suchthat \cont(t) \neq 0}. \nn
\end{align}
Then the \(\lpL[0]\)-norm of \(\cont(\cdot)\) is given by 
\begin{align}
    \|u(\cdot)\|_0 \Let \mu \bigl( \text{supp}(u) \bigr) = \int_0^{\horizon} \abs{\cont(t)}^0 \odif{t},
\end{align}
where the zero exponent is defined by \(z^0 \Let 1\) if \(z\neq 0\) and \(0\) if \(z = 0\).\footnote{Note that \(\norm{\cont(\cdot)}_0\) does not satisfy the absolute homogeneous property and thus, technically, its a semi-norm.} Thus \(\norm{\cont(\cdot)}_0\) is the time length over which the function \(\cont(\cdot)\) admits nonzero values. Over the set of feasible controls, we consider the \highlight{\(\lpL[0]\) robust optimal control problem} which is the central object of this article: 
\begin{myOCP}
    \begin{equation}
	\label{e:OCP} 
\begin{aligned}
& \inf_{\cont(\cdot)}	&& \objlzero \bigl(\cont(\cdot)\bigr) \Let \norm{\cont}_0 = \int_{\tinit}^{\horizon} \abs{\cont(t)}^0 \odif{t} \\
&  \sbjto		&&  \begin{cases}
\text{dynamics}\,\,\eqref{eq:sys},\,\st(\tinit;\uparam)= \param,\\ 
\cont(t) \in \admcont \,\,\text{for a.e}\,\, t\in \lcrc{0}{\horizon},\\
\st(T;\uparam) \in C \text{ for all }\uparam \in \pset.
\end{cases}
\end{aligned}
\end{equation}
\end{myOCP}
This is the robust version of the typical setting of the maximum hands-off control paradigm in the nominal noiseless/uncertainty-free scenario, where the final objective is to obtain a \emph{sparse} optimal control that drives the final state to zero, satisfying all the constraints. Our setting is more challenging: in the presence of uncertainties, and a possible uncountable family of constraints, steering the final state to zero is not practical and thus the OCP \eqref{e:OCP} considers steering \(\st(\horizon;\uparam)\) to a prespecified set for all possible realizations of the uncertainty \(\uparam \in \pset\).

\subsection{Relaxation}
The optimization problem \eqref{e:OCP} is an infinite-dimensional, nonsmooth, and nonconvex robust optimization problem with an uncountable family of constraints, rendering it extremely challenging to solve numerically. 
Following the standard approach of the literature on nominal maximum hands-off control \cite{ref:MH-book-20,ref:MN:MaxHandsOff,ref:MN:23,ref:TI:KK-19}, we relax the problem \eqref{e:OCP} into an \(\lpL[1]\) robust optimal control problem to enable further analysis. This approach also facilitates the use of a broad spectrum of techniques from robust and semi-infinite optimization. With this motivation, consider the relaxed \highlight{\(\lpL[1]\) robust optimal control problem}:
\begin{myOCP}
    \begin{equation}
	\label{e:OCP:L1} 
\begin{aligned}
& \inf_{\cont(\cdot)}	&& \objlone \bigl(\cont(\cdot)\bigr) \Let \norm{\cont(\cdot)}_{\lpL[1]}= \int_{\tinit}^{\horizon} \abs{\cont(t)} \odif{t} \\
&  \sbjto		&&  \begin{cases}
\text{dynamics}\,\,\eqref{eq:sys},\,\st(\tinit;\uparam)= \param,\\ 
\cont(t) \in \admcont \,\,\text{for a.e}\,\, t\in \lcrc{0}{\horizon},\\
\st(T;\uparam) \in C \text{ for all }\uparam \in \pset.
\end{cases}
\end{aligned}
\end{equation}
\end{myOCP}
The optimal control problem \eqref{e:OCP:L1} is an infinite-dimensional robust optimization problem. In our main results in \S\ref{sec:main_result}, we present a series of theoretical results addressing the existence of solutions, necessary optimality conditions, equivalence properties, and an algorithmic framework designed to compute sparse control actions efficiently and with numerical viability.
\section{Main results}\label{sec:main_result}
Our main results are divided into two major parts: (a) \highlight{Theoretical developments} and (b) \highlight{Algorithmic developments}.
\begin{itemize}[leftmargin=*, label={\highlight{\(\circ\)}}]
    \item On the theoretical front, we show the existence of optimal solutions to the OCP \eqref{e:OCP:L1} and consequently we derive a robust (minmax) Pontryagin-type Maximum Principle. We demonstrate that under standard assumptions, the optimal control exhibits a bang-off-bang structure. Subsequently, we prove that the set of optimal solutions for the OCP \eqref{e:OCP:L1} is equivalent to that of the OCP \eqref{e:OCP}, establishing the equivalence of the two problems.
    \item On the algorithmic front, by exploiting the established equivalence and demonstrating that \eqref{e:OCP:L1} is a convex semi-infinite program, we establish an algorithmic framework to solve \eqref{e:OCP:L1} (and, consequently, \eqref{e:OCP}) in an \emph{exact fashion}.
\end{itemize} 
\subsection{Theoretical development: maximum principle and equivalence}\label{subsec:theo:dev:pmp}
We present our main results in the following sequence:
\begin{enumerate}[leftmargin=*, label={\highlight{\(\circ\)}}]
    \item We begin by showing an existence result for the \(\lpL[1]\) robust OCP \eqref{e:OCP:L1}; this is documented in Theorem \ref{thrm:L1:existence}. 
    \item We establish a robust Pontryagin-type maximum principle for the \(\lpL[1]\) robust OCP \eqref{e:OCP:L1}. Our primary tool relies on a nonsmooth robust maximum principle from \cite{ref:minmax:ness:vinter}. This is documented in Theorem \ref{thrm:PMP:conditions}.
    \item Using the structure of the control obtained from Theorem \ref{thrm:PMP:conditions}, we show that under some additional hypothesis on the problem data \(\lpL[1]\) optimal controls admit a bang-off-bang-type representation. This is documented in Corollary \ref{corr:bang-off-bang}.
    \item Finally, leveraging these results, we establish that the OCPs \eqref{e:OCP} and \eqref{e:OCP:L1} are equivalent in the sense that the set of optimal solutions of the \(\lpL[0]\) problem and the \(\lpL[1]\) problem are identical. The result is given in Theorem \ref{thrm:equivalence}.
\end{enumerate}


\begin{figure}[h!]
\centering
\begin{tikzpicture}[
    node distance=1.7cm and 1.6cm,
    block/.style={
        draw,
        rounded corners,
        minimum width=3.5cm,
        text width=3.2cm,
        minimum height=2cm,
        align=center,
        drop shadow,
        fill=#1,
        fill opacity=0.8,
        text opacity=1
    },
    arrow/.style={
        ->,
        line width=1.2pt,
        draw=black!60,
        -Latex
    }
]

\node[block=Orchid!30!white] (L0) {The Robust \\ \(\lpL[0]\) Problem};
\node[block=Bittersweet!30!white, below=of L0] (L1) {The Robust \\ \(\lpL[1]\) Problem};
\node[block=Goldenrod!30!white, below=of L1] (PMP) {Robust PMP \\ (Theorem \ref{thrm:PMP:conditions})};

\node[block=CornflowerBlue!25!white, right=of L0] (Existence) {Existence of \\ \(\lpL[1]\) Solutions \\ (Theorem \ref{thrm:L1:existence})};
\node[block=LimeGreen!25!white, below=of Existence] (BOB) {Bang-Off-Bang \\ Structure \\ (Corollary \ref{corr:bang-off-bang})};
\node[block=Teal!25!white, below=of BOB] (Equivalence) {\(\lpL[0]/\lpL[1]\) Equivalence \\ (Theorem \ref{thrm:equivalence})};

\path[arrow] (L0) edge node[midway, above, fill=white, inner sep=1pt] {Relaxation} (L1);
\path[arrow] (L1) edge (PMP);
\path[arrow] (L1) edge (Existence);
\path[arrow] (PMP) edge (BOB);
\path[arrow] (BOB) edge (Equivalence);
\path[arrow, draw=black!50, opacity=0.7] (Existence) edge[bend right=20] (Equivalence);

\end{tikzpicture}
\caption{Roadmap of the robust hands-off synthesis. The original OCP~\eqref{e:OCP}, featuring an \(\lpL[0]\) objective and semi-infinite uncertainty constraints, is reduced --- via the \(\lpL[0]\)–\(\lpL[1]\) equivalence (Theorem~\ref{thrm:PMP:conditions} and Theorem~\ref{thrm:equivalence}, under the conditions of Corollary~\ref{corr:bang-off-bang}) --- to a convex robust \(\lpL[1]\) formulation. For practical implementation, in \S \ref{subsec:algo:dev:sip}, after performing piecewise constant control parametrization, the problem is further transformed into a finite-dimensional convex semi-infinite program~\eqref{eq:SR:pre_SIP}, which is then addressed using tools from robust optimization.
}
\label{fig:log:flow}
\end{figure}
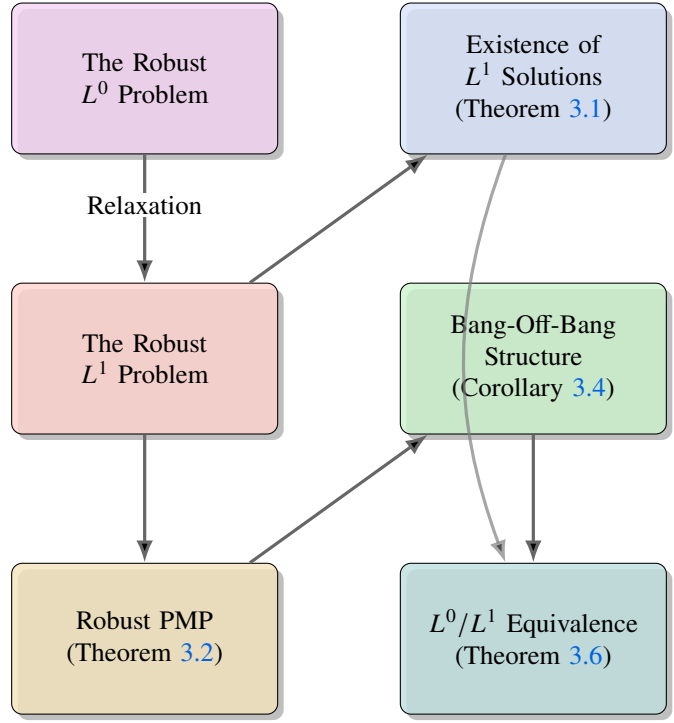

We begin with the existence theorem. 
\begin{myOCP}
\begin{theorem}\label{thrm:L1:existence}
Consider the optimal control problems \eqref{e:OCP}--\eqref{e:OCP:L1} along with their data \ref{ocp:data:1}--\ref{ocp:data:3}. Define the set of admissible control actions
\begin{align}\label{eq:adm_cont_set}
    \admu \Let \left\{\cont(\cdot) \in \lpL[\infty](\lcrc{0}{\horizon};\Rbb) \;\middle\vert\;  
    \begin{array}{@{}l@{}}
        \cont(t) \in \admcont \text{ for a.e. }t \in \lcrc{0}{\horizon} \text{ and }\\
        \psi\bigl(\st(\horizon;\uparam)\bigr) \le 0 \text{ for all }\uparam \in \pset
        \end{array}
        \right\}. 
\end{align}
Note that for a given \(\uparam \in \pset\) and \(\st(0;\uparam) = \param\), by \(\st_{u}(\horizon;\uparam)\) we denote the state trajectory of the ODE \eqref{eq:sys}, under the control action \(\cont(\cdot)\), evaluated at time \(\horizon\).\footnote{The state trajectory at \(t = \horizon\) is given by (employing Duhamel's formula)
\begin{align}
 \st_{\cont}(\horizon,\uparam)\Let \epower{A(\uparam)\horizon}\param + \int_{0}^{\horizon} \epower{A(\uparam)(s-\horizon)}b(\uparam)\cont(s)\odif{s} \teL \mathsf{E}_{\param,\horizon,\uparam}(\cont(\cdot)),
\end{align}
where the quantity \(\mathsf{E}_{\param,\horizon,\uparam}(\cont(\cdot))\) is the \emph{end-point mapping}; see \cite[Definition 1.1]{ref:trelat2023control}.} Let \(\admu\) be nonempty, then the OCP \eqref{e:OCP:L1} admits a solution. 
\end{theorem}
\end{myOCP}
\begin{proof}
We provide a proof using the direct method in the calculus of variations \cite[Chapter 3]{ref:santambrogio2023course}. Define 
\begin{align}
   \objective\as \Let \inf_{\cont(\cdot) \in \admu} \objlone \bigl(\cont(\cdot)\bigr),
\end{align}
and note that since the admissible set \(\admu\) is nonempty \(\objective\as < +\infty.\) Then, by definition of the infimum there exists a \emph{minimizing sequence} \(\bigl(\cont_k(\cdot)\bigr)_{k \in \N} \subset \admu\) such that \(\lim_{k \lra +\infty} \objlone \bigl(\cont_k(\cdot)\bigr) = \objective\as\). Moreover, by definition
\[\bigl(\cont_k(\cdot)\bigr)_{k \in \N} \subset \Ball_{\infty} \Let \aset[\big]{\mu(\cdot) \in \lpL[\infty]\bigl(\lcrc{0}{\horizon};\Rbb\bigr) \suchthat \norm{\mu(\cdot)}_{\infty} \le 1}.\]
The celebrated Banach-Alaoglu theorem \cite[Theorem 3.14]{ref:clarke_ocpbook} states that the closed unit ball of the dual of a given Banach space is sequentially compact in the weak\(\as\) topology \cite[\S 3.3]{ref:clarke_ocpbook}. Since \(\lpL[\infty]\bigl(\lcrc{0}{\horizon};\Rbb\bigr)\) is the dual of the separable space \(\lpL[1]\bigl(\lcrc{0}{\horizon};\Rbb\bigr)\), this theorem applies. Therefore, there exists a subsequence (which we do not relabel for simplicity) and a mapping \(\widehat{\cont}(\cdot) \in \lpL[\infty]\bigl(\lcrc{0}{\horizon};\Rbb\bigr)\) with \(\norm{\widehat{\cont}(\cdot)}_{\infty} \le 1\) such that\footnote{
In our setting, admissible controls satisfy
\[
u(\cdot)\in \lpL[\infty]([0,T];\mathbb{R}) \quad \text{with }\, u(t) \in \mathbb{U} \text{ for a.e. }t \in [0,T].
\]
The key functional-analytic facts we use are:
\begin{itemize}[label = \(\circ\),leftmargin=*]
\item $\lpL[\infty]([0,T];\mathbb{R})$ is (isometrically) the dual of the separable Banach space $\lpL[1]([0,T];\mathbb{R})$ \cite[Chapter 5 and Chapter 6, Theorem 6.10]{ref:clarke_ocpbook};
\item the closed unit ball in $\lpL[\infty]([0,T];\Rbb)$ is therefore sequentially weak$^*$-compact by the Banach--Alaoglu theorem \cite[Chapter 3, Theorem 3.14 and Chapter 6]{ref:clarke_ocpbook}.
\end{itemize}
We provide the definition of weak$^*$ convergence in $\lpL[\infty]([0,T];\Rbb)$ \cite[Chapter 3]{ref:clarke_ocpbook}): A sequence $\bigl(u_k(\cdot)\bigr)_{k \in \N} \subset \lpL[\infty]([0,T];\Rbb)$ converges to $\bar u(\cdot)\in\lpL[\infty]([0,T];\Rbb)$ in the \emph{weak$^*$-sense} if
\[
\int_0^T u_k(t) f(t)\,dt \;\longrightarrow\; \int_0^T \bar u(t) f(t)\,dt
\quad \text{for every } f(\cdot)\in L^1([0,T];\Rbb) \text{ as } k \lra +\infty.
\]
} 
\[
\cont_k(\cdot) \xrightarrow{\text{weak}^*} \widehat{u}(\cdot) \quad\text{in } \lpL[\infty]\bigl(\lcrc{0}{\horizon};\Rbb\bigr),
\]
which implies that for any \(f(\cdot) \in \lpL[1]\bigl(\lcrc{0}{\horizon};\Rbb\bigr)\), we have
\begin{align}
\int_0^{\horizon} \cont_k(t) f(t) \odif{t} \xrightarrow{k \lra + \infty}\int_0^{\horizon} \widehat{\cont}(t) f(t) \odif{t}.  
\end{align}
We show convergence of the state trajectories. Fix \(\uparam \in \pset\) and define
\begin{align}
  f_k(\uparam) \Let \st_{\cont_k}(\horizon;\uparam) = \epower{A(\uparam) \horizon} \param + \int_0^{\horizon} \epower{A(\uparam)(\horizon - s)} b(\uparam) \cont_k(s)\odif{s}.  
\end{align}
Pointwise convergence for fixed \(\uparam \in \pset\) follows easily; indeed: the function \(s \mapsto \epower{A(\uparam)(\horizon - s)} b(\uparam)\) is a continuous function defined on a compact interval \(\lcrc{0}{\horizon}\) and thus it is integrable. By the definition of weak\(\as\) convergence we see that
\begin{align}
\lim_{k \lra +\infty} \int_0^{\horizon} \epower{A(\uparam)(\horizon - s)} b(\uparam) \cont_k(s)\odif{s} = \int_0^{\horizon} \epower{A(\uparam)(\horizon - s)} b(\uparam) \widehat{u}(s)\odif{s}.
\end{align}
Thus, for fixed \(\uparam \in \pset\) we have
\begin{align}\label{eq:pwise:weak-star}
f_k(\uparam) = x_{u_k}(T;\uparam) \xrightarrow{k \lra +\infty} x_{\widehat{u}}(T;\uparam) \teL \widehat{f}(\uparam).
\end{align}
To show this convergence is \emph{uniform} on \(\pset\), we need to invoke the Arzela-Ascoli theorem, and to this end, we check if the family \((f_k)_{k\in \N}\) is uniformly bounded and equicontinuous. Uniform boundedness follows immediately because \(\pset \ni \uparam \mapsto A(\uparam)\) and \(\pset \ni \uparam \mapsto b(\uparam)\) are continuous on the compact set \(\pset\) and \(\norm{\cont_k(\cdot)}_{\infty} \le 1\). To show equicontinuity, we show \cite[Chapter 11, Definition 11.27]{ref:RudRealCom-1987} that for all \(\eps>0\), there exists \(\delta>0\), for any \(\uparam_1,\uparam_2 \in \pset\) with \(\norm{\uparam_1 -\uparam_2} < \eps\), \(\norm{f_k(\uparam_1) - f_k(\uparam_2)} < \eps\). Observe that
\begin{align}
    \norm{f_k(\uparam_1) - f_k(\uparam_2)} &\le \norm{\epower{A(\uparam_1)\horizon} - \epower{A(\uparam_2)\horizon}} \norm{\param} \nn \\& + \int_0^{\horizon} \norm{\epower{A(\uparam_1)(\horizon-s)}b(\uparam_1) - \epower{A(\uparam_2)(\horizon-s)}b(\uparam_2)}\odif{s}
\end{align}
for \(s \in \lcrc{0}{\horizon}\). Since \(\uparam \mapsto \epower{A(\uparam)\horizon}\) and \(\uparam \mapsto \epower{A(\uparam)(\horizon-s)}b(\uparam)\) are uniformly continuous on \(\pset\), for any \(0<\eps'\Let \frac{\eps}{\norm{\param}+\horizon}\), there exists a \(\delta>0\) such that with \(\norm{\uparam_1 - \uparam_2} < \delta\) we have \(\norm{\epower{A(\uparam_1)\horizon} - \epower{A(\uparam_2)\horizon}} < \eps'\) and \(\norm{\epower{A(\uparam_1)(\horizon-s)}b(\uparam_1) - \epower{A(\uparam_2)(\horizon-s)}b(\uparam_2)} < \eps'\) for all \(s \in \lcrc{0}{\horizon}\). Consequently, 
\begin{align}
    \norm{f_k(\uparam_1) - f_k(\uparam_2)} \le \eps'\bigl( \norm{\param} + \horizon \bigr) = \eps. 
\end{align}
Thus, the family \((f_k)_{k\in \N}\) is uniformly bounded and equicontinuous on \(\pset\); and consequently, by Arzela-Ascoli theorem \cite[Chapter 11, Theorem 11.28]{ref:RudRealCom-1987} there exists a subsequence (not relabeled) such that \(f_k \to \widehat f\) uniformly on \(\pset\), where the limit \(\widehat f(\uparam)=x_{\widehat u}(\horizon;\uparam)\) is the pointwise limit given by \eqref{eq:pwise:weak-star}.

We are now ready to show that \(\widehat{\cont}(\cdot) \in \admu\). By definition of the admissible set in \eqref{eq:adm_cont_set}, for every \(\cont_k(\cdot) \in \admu\) for \(k \in \N\), we have \(\psi \bigl(\st_{\cont_k}(\horizon;\uparam) \bigr) \le 0\) for all \(\uparam \in \pset\). We just demonstrated the uniform convergence 
\[
\st_{\cont_k}\bigl(\horizon;\uparam\bigr) \xrightarrow{\norm{\cdot}_{\infty}} \st_{\widehat{\cont}}\bigl(\horizon;\uparam\bigr)\quad \text{ as }k\lra +\infty,
\]
and by design \(\psi(\cdot)\) is continuous. Passing to the limit, we get 
\begin{align}
\lim_{k \lra +\infty} \psi\bigl( \st_{\cont_k}(\horizon; \uparam)\bigr) = \psi\biggl(\lim_{k \lra +\infty} \st_{\cont_k}(\horizon; \uparam)\biggr) = \psi\bigl(\st_{\widehat{\cont}}(\horizon; \uparam)\bigr),   
\end{align}
and \(\psi\bigl(\st_{\widehat{\cont}}(\horizon;\uparam)\bigr) \leq 0\) for all \(\uparam \in \pset\). Thus \(\widehat{\cont}(\cdot) \in \admu\).

Finally, we show that \(\widehat{\cont}(\cdot)\) is an optimal choice. To this end, note that \(\cont(\cdot) \mapsto \objlone\bigl(\cont(\cdot)\bigr) \Let  \int_0^{\horizon} \abs{\cont(t)} \odif{t}\) is convex and continuous and thus it is weak\(\as\) lower semicontinuous \cite[Theorem 2.12]{ref:FT:PDE:OptConBook}. Then
\begin{align}\label{ef:objective:ineq:1}
\objlone\bigl(\widehat{\cont}(\cdot)) \le \liminf_{k \lra + \infty} \objlone\bigl( \cont_k(\cdot) \bigr) = \objective\as.
\end{align}
However, \(\widehat{\cont}(\cdot) \in \admu\) is a feasible control and thus \(\objlone \bigl( \widehat{\cont}(\cdot) \bigr) \ge \objective\as\) which implies, together with \eqref{ef:objective:ineq:1}, that \(\objlone \bigl( \widehat{\cont}(\cdot) \bigr) = \objective\as\). Thus, \(\widehat{\cont}(\cdot)\) is an \(\lpL[1]\) optimal control and the proof is complete. 
\end{proof}

We now establish a robust maximum principle for the OCP \eqref{e:OCP:L1}.
\begin{myOCP}
\begin{theorem}\label{thrm:PMP:conditions}
Consider the optimal control problem \eqref{e:OCP:L1} along with its data \ref{ocp:data:1}--\ref{ocp:data:3}. Let the optimal trajectory-pair \(\bigl(\cont\as(\cdot),\aset[]{\st\as(\cdot;\uparam) \suchthat \uparam \in \pset}\bigr)\) be a strong local minimizer of \eqref{e:OCP:L1} and let the hypothesis of Theorem \ref{thrm:L1:existence} holds. Then there exists:
\vspace{-2mm}
\begin{enumerate}[leftmargin=*, widest=b, align=left]
    \item A scalar \(\multip \in \lcrc{0}{1}\);
    \item A non-negative Radon measure \(\measure\) on the compact set \(\pset\);
    \item A family of absolutely continuous functions called the costate or adjoint arcs \(\adjx(\cdot; \uparam) \in W^{1,1}\bigl(\lcrc{0}{\horizon};\Rbb^{d}\bigr)\), defined for \(\measure\)-a.e. \(\uparam \in \supp(\measure)\),
\end{enumerate}
\vspace{-1mm}
satisfying the following conditions:
\begin{enumerate}[label=\textup{{\highlight{(\ref{thrm:PMP:conditions}-\alph*)}}}, leftmargin=*, widest=b, align=left]
    \item \label{pmp:cond:1} (Nontriviality) The multipliers are not simultaneously zero and are normalized such that:
   \begin{align}
     \multip + \measure(\pset) = 1.  \nn 
   \end{align}

    \item \label{pmp:cond:2} (Adjoint dynamics) For \(\measure\)-a.e. \(\uparam \in \pset\), the costate arc \(\adjx(\cdot; \uparam)\) is a solution to the adjoint equation:
    \begin{align}
    -\dot{\adjx}(t; \uparam)  = A(\uparam)^{\top} p(t; \uparam) \quad \text{for a.e. } t \in \lcrc{0}{\horizon}. \nn     
    \end{align}

    \item \label{pmp:cond:3} (Transversality condition) For \(\measure\)-a.e. \(\uparam \in \pset\), there exists a scalar \(\mu(\uparam) \ge 0\) such that the final value of the costate is given by:
    \begin{align}
        -\adjx(\horizon; \uparam) = \mu(\uparam) \nabla\psi\bigl(\st\as(\horizon; \uparam)\bigr). \nn
    \end{align}

     \item \label{pmp:cond:4} (Support and complementarity condition) The measure \(\measure\) is supported on the set of "worst-case" parameters for which the terminal constraint is active:
    \begin{align}
        \supp(\measure) \subset \aset[\big]{\uparam \in \pset \suchthat \psi\bigl(\st\as(\horizon; \uparam)\bigr) = 0}. \nn 
    \end{align}
    Moreover, for every \(\uparam \in \supp(\measure)\), \(\mu(\uparam)>0\) only if \(\psi\bigl(\st\as(\horizon;\uparam)\bigr) = 0 \). 

    \item \label{pmp:cond:5} (Hamiltonian maximization condition) For almost every \(t \in \lcrc{0}{\horizon}\), the optimal control \(\cont\as(\cdot)\) is a solution to the maximization problem
    \begin{align}
      \cont\as(t) & \in \argmax_{v \in \admcont} \biggl(\int_{\pset}\adjx(t;\uparam)^{\top}\bigl(A(\uparam)\st\as(t;\uparam)+b(\uparam)v\bigr)\measure(\odif{\uparam}) - \multip \abs{v} \biggr)\nn  \\& =\argmax_{v \in \admcont} \aset[\big]{\switch(t) v - \eta \abs{v}}  \nn
    \end{align}
    where the averaged switching function \(t \mapsto \switch(t) \in \Rbb\) is defined by
    \[ \lcrc{0}{\horizon} \ni t \mapsto \switch(t) \Let \int_{\pset} p(t; \uparam)^{\top} b(\uparam) \measure(\odif{\uparam}) \in \Rbb.\]
\end{enumerate}
\end{theorem}
\end{myOCP}

The proof proceeds through four key steps. First, the \(\lpL[1]\) OCP \eqref{e:OCP:L1} is reformulated to the standard Mayer form. Next, by introducing an artificial parameter, an equivalent minmax formulation is derived. This formulation leverages results from \cite{ref:minmax:ness:vinter} to establish a set of necessary conditions. Finally, these conditions are refined to obtain the necessary conditions for the OCP given in \eqref{e:OCP:L1}.

\begin{proof}
We reformulate the OCP \eqref{e:OCP:L1}, which is in the standard Bolza form, to the Mayer form. Towards this end, we augment the trajectory $t \mapsto x(t) \in \Rbb^d$ with $t \mapsto x_{d+1}(t) \in \Rbb$ where
\begin{align}
\dot{x}_{d+1}(t) = \abs{u(t)}, \quad x_{d+1}(0) = 0 \text{ for all }t \in \lcrc{0}{\horizon}. \nn
\end{align} 
We will denote the augmented state trajectory by \(t \mapsto y(t) \Let [x(t)^{\top}, x_{d+1}(t)]^{\top} \in \Rbb^{d+1}\). Thus OCP \eqref{e:OCP:L1} is now equivalent to an OCP where the objective function consists of only the terminal cost \(x_{d+1}(\horizon)\) with the augmented dynamics given by 
\begin{align}\label{eq:aug:dyn}
\dot{y}(t;\uparam) = F(y(t;\uparam),u(t)) \Let \begin{pmatrix} A(\uparam)\st(t) + b(\uparam)\cont(t) \\ \abs{u(t)} \end{pmatrix}, \quad y(0;\uparam) \Let \begin{pmatrix}
    \param \\ 0
\end{pmatrix}. 
\end{align}
The OCP \eqref{e:OCP:L1} becomes
\begin{equation}\label{e:OCP:L1:Mayer} 
\begin{aligned}
& \inf_{\cont(\cdot)}	&& g\bigl(y(\horizon;\uparam),\uparam\bigr) \Let x_{d+1}(\horizon) \\
&  \sbjto		&&  \begin{cases}
\text{dynamics}\,\,\eqref{eq:aug:dyn},\,y(\tinit;\uparam)= (\param,0)^{\top},\\ 
\cont(t) \in \admcont \,\,\text{for a.e}\,\, t\in \lcrc{0}{\horizon},\\
y(T;\uparam) \in C \Let \aset[\big]{(\xi_1,\xi_2) \in \Rbb^d \times \Rbb\suchthat \psi(\xi_1) \le 0} \text{ for all }\uparam \in \pset.
\end{cases}
\end{aligned}
\end{equation}
We recast the Mayer-form OCP \eqref{e:OCP:L1:Mayer} into an equivalent minmax problem by extending the uncertainty set via augmenting an artificial parameter \(\uparam\as\). Define \(\apset \Let \pset \cup \aset[]{\uparam\as}\); consider the objective function
\(g(\augst,\widetilde{\uparam})\) defined by
\begin{align}
g(\augst, \widetilde{\uparam}) \Let \begin{cases} \st_{d+1} & \text{if } \widetilde{\uparam} = \uparam\as \\ 0 & \text{if } \widetilde{\uparam} \in \pset, \end{cases} \nn
\end{align}
and a new terminal constraint, \(C(\widetilde{\uparam})\) by
\begin{align}
C(\widetilde{\uparam}) \Let \begin{cases} \aset[\big]{\augst = (\st,\st_{d+1}) \in \Rbb^{d+1} \suchthat \psi(\st) \le 0} & \text{if } \widetilde{\uparam} \in \pset \\ \Rbb^{d+1} & \text{if } \widetilde{\uparam} = \uparam\as. \end{cases}    \nn 
\end{align}
Thus, the equivalent minmax reformulation of \eqref{e:OCP:L1:Mayer} is
\begin{equation}\label{eq:reformulated:minmax} 
\begin{aligned}
& \inf_{\cont(\cdot)\in \admcon} \sup_{\widetilde{\uparam} \in \apset} && g\bigl (\augst(\horizon;\widetilde{\uparam}), \widetilde{\uparam}\bigr)  \\
&  \sbjto		&&  \begin{cases}
\augst(\horizon;\widetilde{\uparam}) \in C(\widetilde{\uparam}) \text{ for all }\widetilde{\uparam} \in \apset.
\end{cases}
\end{aligned}
\end{equation}
Define \(\pset_1 \Let \pset\) and \(\pset_2 \Let \uparam\as\); and thus \(\apset = \pset_1 \cup \pset_2\). Our next step involves application of the robust PMP in \cite{ref:minmax:ness:vinter}. To this end, first, we need to verify the conditions \cite[S1-S5, pp.946-947]{ref:minmax:ness:vinter}. Observe that 
\begin{itemize}[leftmargin=*]
    \item The vector field \((\augst,\cont,\uparam) \mapsto F(\augst, \cont)\) is measurable for each \(\augst \in \Rbb^{d+1}\), Thus, (S1) holds. 

    \item the Jacobian of \(F(\cdot)\) with respect to \(\augst\) is \(\nabla_{\augst} F =  \begin{psmallmatrix} A(\uparam) & 0 \\ 0 & 0 \end{psmallmatrix}\) and recall that \(\uparam \mapsto A(\uparam)\) is continuous on the compact set \(\pset\). Thus \(\augst \mapsto F(\augst,\cont)\) is Lipschitz for all \(\uparam \in \pset\). Furthermore, since all problem data and the control set \(\admcont\) are bounded, \(F(\cdot)\) is also bounded, and (S2) holds. 

    \item  Recall that, for any \(\uparam \in \pset\) the function \(g(\augst, \uparam)\) is identically zero, and thus (S3) is satisfied. 

    \item (S4) holds because \(\pset \ni \uparam \mapsto A(\uparam)\) and \(\pset \ni \uparam \mapsto b(\uparam)\) are continuous and (S5) holds trivially. 
\end{itemize}
An application of \cite[Theorem 3.2]{ref:minmax:ness:vinter} guarantees the existence of a Radon probability measure \(\widetilde{\measure}\) on the set \(\apset\), and a family of arcs 
\begin{align}\label{eq:costate:arcs}
    \adj(\cdot;\widetilde{\uparam}) \Let \bigl( \adjx(\cdot;\widetilde{\uparam}),\adjy(\cdot;\widetilde{\uparam}\bigr) \in W^{1,1}\bigl( \lcrc{0}{\horizon};\Rbb^{d+1}\bigr)\text{ for } \widetilde{\measure}\text{-a.e. } \widetilde{\uparam}.
\end{align}
For the OCP \eqref{eq:reformulated:minmax} let the pair \(\bigl(\augst\as(\cdot;\uparam),\cont(\cdot)\bigr)\) denote the optimal state-action trajectory. Define the Hamiltonian
\begin{align}
\Rbb^{d+1} \times \Rbb^{d+1} \times \Rbb \ni (\augst,\adj,\cont) & \mapsto \aughamil(\augst,\adj,\cont) \Let \inprod{\adj}{F(\augst,\cont)} \nn \\& = \inprod{\adjx}{A(\uparam)\st + b(\uparam)\cont} + \adjy\abs{\cont} \in \Rbb, \nn
\end{align}
where \(t\mapsto \adj(t;\uparam) \Let (\adjx(t;\uparam)^{\top}, \adjy(t;\uparam))^{\top} \in \Rbb^{d+1}\) for all \(t \in \lcrc{0}{\horizon}\) is the augmented costate/adjoint defined in \eqref{eq:costate:arcs}. Now we will distil the conditions \ref{pmp:cond:2}--\ref{pmp:cond:5} from the necessary conditions for the OCP \eqref{eq:reformulated:minmax}.

We have, by definition, \(\widetilde{\measure}(\apset) = 1\); we split \(\widetilde{\measure}\) into its restriction on the set \(\pset\) and on the point mass \(\aset[]{\uparam\as}\). We define our multipliers as: \(\eta \Let \tilde{\measure}(\aset[]{\uparam\as})\) (the mass of the measure on the artificial point) and \(\measure\) is the restriction of \(\tilde{\measure}\) to the set \(\pset\). Since \(\apset\) is a disjoint union of \(\pset\) and \(\aset[]{\uparam\as}\), we can split the measure's mass and write \(\tilde{\measure}(\pset) + \tilde{\measure}(\aset[]{\uparam\as}) = 1\), giving us the nontriviality condition \(\measure(\pset) + \eta = 1\).

Let us derive the condition \ref{pmp:cond:2}--\ref{pmp:cond:4}. For the OCP \eqref{eq:reformulated:minmax}, corresponding to the pair \(\bigl(\augst\as(\cdot;\uparam),\cont(\cdot)\bigr)\), and the function
\(\aughamil(\cdot)\), we have the adjoint equation: \(-\dot{\adj} = \nabla_{\augst} \aughamil\), from which, computing the gradients
\begin{align}
\nabla_{\augst} \aughamil = \begin{pmatrix} \nabla_{\st} \aughamil \\ \nabla_{\st_{d+1}} \aughamil \end{pmatrix} = \begin{pmatrix} A(\uparam)^{\top} \adjx \\ 0 \end{pmatrix}, \nn  
\end{align}
we immediately obtain the adjoint equations 
\begin{align}\label{eq:adj:1:2}
-\dot{\adjx}(t; \uparam) &= A(\uparam)^{\top} \adjx(t; \uparam), \nn \\ -\dot{p}_{d+1}(t; \uparam) &= 0, \text{ which implies } \adjy(t;\uparam) \text{ is constant in time.}
\end{align}

The support and the transversality conditions follow directly from the more general conditions in \cite[Theorem 3.2]{ref:minmax:ness:vinter}. We have 
\begin{align}\label{eq:gen:transversality}
-\adj(\horizon;\widetilde{\uparam}) \in N_{C(\widetilde{\uparam})}\bigl(\augst\as(\horizon;\widetilde{\uparam})\bigr) + \partial_{\augst} g\bigl(\augst\as(\horizon;\widetilde{\uparam}), \widetilde{\uparam}\bigr).  
\end{align}
Observe that, when \(\uparam \in \pset\), \(g(y, \uparam) = 0\), and thus \eqref{eq:gen:transversality} reduces to
\begin{align}
-\adj(\horizon;\uparam) \in N_{\aset[]{\augst \suchthat \psi(x) \le 0}}\bigl(\augst\as(\horizon;\uparam)\bigr)    
\end{align}
The normal cone to this sublevel set is non-zero only if the point is on the boundary, i.e., \(\psi\bigl(\st\as(\horizon;\uparam)\bigr)=0\). Since the measure \(\measure\) can only be supported where the conditions are active, this directly proves the condition \ref{pmp:cond:4}.

From the adjoint equations and the preceding arguments, it follows that there exists a scalar multiplier \(\mu(\uparam) \ge 0\) such that
\begin{align}
-\adj(\horizon;\uparam) = -\begin{pmatrix} \adjx(\horizon;\uparam) \\ \adjy(\horizon;\uparam) \end{pmatrix} = \mu(\uparam) \begin{pmatrix} \nabla_{\st}\psi\bigl(\st\as(\horizon;\uparam)\bigr) \\ 0 \end{pmatrix},\nn
\end{align}
where the first equation is the transversality condition in \eqref{pmp:cond:3} and the second equation implies that \(\adjy(\horizon;\uparam) = 0\), and this implies. from \eqref{eq:adj:1:2}, that \(\adjy(t;\uparam)\) is identically zero for all \(t \in \lcrc{0}{\horizon}\), when \(\uparam \in \pset\). 

When \(\widetilde{\uparam} = \uparam\as\), the constraint set is \(C(\uparam\as) = \Rbb^{d+1}\), and thus \(N_{C(\uparam\as)} = \aset[]{0}\). The cost is \(g(\augst, \uparam\as) = \st_{d+1}\) and thus \eqref{eq:gen:transversality} reduces to
\begin{align}
-\adj(\horizon;\uparam\as) \in \aset[]{0} + \aset[\big]{[0, \ldots, 0, 1]^{\top}}
\end{align}
which implies \(-\adjx(\horizon;\uparam\as) = 0\) and \(-\adjy(\horizon;\uparam\as) = 1\). From the adjoint dynamics \eqref{eq:adj:1:2} we have \(\adjx(t;\uparam\as) = 0\) and \(\adjy(t;\uparam) = -1\) for all \(t \in \lcrc{0}{\horizon}\). 

Finally, the minmax PMP states that the control trajectory \(\cont\as(\cdot)\) must solve the variational problem 
\[\max_{\cont \in \admcont} \int_{\apset} \widehat{H}_{\widetilde{\uparam}} \bigl( \augst\as(t;\widetilde{\uparam}),\adj(t;\widetilde{\uparam}),\cont \bigr) \widetilde{\measure}(\odif{\widetilde{\uparam}}).\]
Expanding upon the integral, we obtain
\begin{align}
&\max_{v \in \admcont} \int_{\pset}\widehat{H}_{{\uparam}} \bigl( \augst\as(t;{\uparam}),\adj(t;{\uparam}),\cont \bigr) \measure(\odif{\uparam})+ \multip \widehat{H}_{\uparam\as}\bigl( \augst\as(t;{\uparam}\as),\adj(t;{\uparam}\as),\cont \bigr) \nn \\
&  \stackrel{\mathclap{(\dag)}}{=} \max_{v \in \admcont} \int_{\pset} \bigl( \adjx(t;\uparam) (A(\uparam)\st\as(t;\uparam) + b(\uparam)v) + \adjy(t;\uparam)\abs{v} \bigr) \measure(\odif{\uparam}) +  \nn \\ & \hspace{20mm}\multip \bigl( \adjx(t;\uparam\as)^{\top} (A(\uparam\as)\st(t;\uparam\as)+b(\uparam\as)v) + \adjy(t;\uparam\as)\abs{v}\bigr) \nn \\
& \stackrel{\mathclap{(\ddag)}}{=} \max_{v \in \admcont} \bigg{\{}\int_{\pset} \bigl( \adjx(t;\uparam)^{\top}(A(\uparam)\st\as(t;\uparam)+b(\uparam)v) + 0 \times \abs{v} \bigr) \measure(\odif{\uparam}) + \nn \\& \hspace{20mm} \multip \big( 0 \times (A(\uparam\as)\st\as(t;\uparam\as)+b(\uparam\as)v) + (-1)\abs{v} \bigr) \bigg{\}}, \nn
\end{align}
where: 
\begin{itemize}[leftmargin=*, label= {\highlight{\(\circ\)}}]
\item In step \((\dag)\): we split the integral over \(\apset\) into integrals over \(\pset\) and \(\aset[]{\uparam\as}\). Recall that \(\multip \Let \widetilde{\measure}(\aset[]{\uparam\as})\) and \(\measure\) is the restriction of the measure \(\widetilde{\measure}\) to \(\pset\).

\item In step \((\ddag)\): recall that, when \(\uparam \in \pset\), we showed that \(\adjy(t;\uparam) \equiv 0\) and for \(\uparam\as\) we showed that \(\adjx(t;\uparam\as) \equiv 0\) and \(\adjy(t;\uparam\as) \equiv -1\) for all \(t \in \lcrc{0}{\horizon}\). 
\end{itemize}
With all these ingredients, the variational problem reduces to
\begin{align}\label{eq:hamimax:prefinal}
\max_{v \in \admcont} \aset[\bigg]{ \int_{\pset} \adjx(t;\uparam)^{\top}(A(\uparam){\st}\as(t;\uparam) + b(\uparam)v) \measure(\odif{\uparam}) - \eta \abs{v}}. 
\end{align}
Since the first term inside the integral in \eqref{eq:hamimax:prefinal} does not depend on \(v \in \admcont\), it does not play any role in the maximization, as thus, defining the averaged switching function \(t \mapsto \switch(t) \in \Rbb\)
    \[ \lcrc{0}{\horizon} \ni t \mapsto \switch(t) \Let \int_{\pset} p(t; \uparam)^{\top} b(\uparam) \measure(\odif{\uparam}) \in \Rbb,\]
we obtain
\begin{align}
\cont\as(t) \in \argmax_{v \in \admcont} \aset[\bigg]{ \left( \int_{\pset} \adjx(t;\uparam)^{\top} b(\uparam) \measure(\odif{\uparam}) \right) v - \eta \abs{v}} = \argmax_{v \in \admcont} \aset[\big]{\switch(t)v - \eta \abs{v}}.
\end{align}
The proof is complete.  
\end{proof}

\begin{remark}
Theorem \ref{thrm:PMP:conditions} is a robust and nonsmooth analogue of Pontryagin’s maximum principle for OCP \eqref{e:OCP:L1} with an \emph{uncountable family of constraints indexed by \(\uparam \in \pset\)}. Under the hypotheses \ref{ocp:data:1}--\ref{ocp:data:3}, the robust PMP yields the existence of (a) a scalar multiplier \(\multip\ge 0\),
(b) a nonnegative Radon measure \(\measure\) (which can be interpreted as a continuous analogue of KKT multipliers over scenarios) supported on the \emph{active uncertainty realizations}, and
(c) a family of adjoint variables \(\adjx(\cdot;\uparam)\) defined \(\measure\)-a.e., such that the optimal control maximizes an averaged Hamiltonian. Finally, the averaged switching function \(\switch(\cdot)\) represents an aggregate sensitivity of the robust terminal feasibility margin with respect to the control input at time \(t\).
\end{remark}




Under some additional assumptions on the problem data, we have the immediate consequence that the \(\lpL[1]\) optimal control is bang-off-bang-type. 
\begin{myOCP}
\begin{corollary}\label{corr:bang-off-bang}
Consider the optimal control problem \eqref{e:OCP}--\eqref{e:OCP:L1} along with their data \ref{ocp:data:1}--\ref{ocp:data:3} and let the hypothesis of Theorem \ref{thrm:L1:existence} hold. We define the \emph{active set} by \
\begin{align}
    \actset \Let \aset[]{ \uparam\in \pset \suchthat \psi\bigl(\st\as(\horizon;\uparam)\bigr)=0}.
\end{align}
denote \(\actmeas \Let \measure|_{\actset}\). Suppose that
\begin{enumerate}[label={\highlight{\textup{(\ref{thrm:PMP:conditions}-\alph*)}}}, leftmargin=*, widest=b, align=left]
\item \label{corr:normality} (Normality condition) The problem \eqref{e:OCP}--\eqref{e:OCP:L1} is normal, i.e., \(\multip>0\).

\item \label{corr:ctrb} (Average controllability-like condition) For every measurable \(\uparam\mapsto v(\uparam)\in \Rbb^d\) with \(v(\uparam) \in N_{C}\bigl(\st\as(\horizon;\uparam)\bigr)\) that is not a.e. zero with respect to \(\actmeas\), there exists \(m \in \aset[]{1,\ldots,d-1}\) such that 
\begin{align}
\int_{\actset} \inprod{A(\uparam)^m b(\uparam)}{v(\uparam)} \actmeas(\odif{\uparam}) \neq\ 0. \nn
\end{align}

\end{enumerate}
Then the optimal control \(t \mapsto \cont\as(t)\) admits a bang-off-bang representation for a.e. \(t \in \lcrc{0}{\horizon}\) and it is given by
\begin{align}
    \lcrc{0}{\horizon} \ni t \mapsto \cont\as(t) \Let \begin{cases} 1 & \text{if } \switch(t) > \multip \\ 0 & \text{if } \abs{\switch(t)} < \multip\\ -1 & \text{if } \switch(t) < -\multip. \end{cases}    \nn 
\end{align}
\end{corollary}
\end{myOCP}

\begin{remark}
Note that, from the transversality condition in Theorem \ref{thrm:PMP:conditions} we have \(-\adjx(\horizon;\uparam)\in N_{C}\bigl(\st\as(\horizon;\uparam)\bigr).\) In particular, if $\uparam$ is \emph{inactive}, i.e., \(\psi\bigl(\st\as(\horizon;\uparam)\bigr) < 0\), then $N_{C(\uparam)}(x^*(T;\uparam))=\{0\}$ and hence \(\adjx(\horizon;\uparam)=0\) for all inactive \(\uparam\). Thus, without loss of generality, we specialized Corollary \ref{corr:bang-off-bang} on \(\actset\). The normality assumption is standard \cite{ref:MN:MaxHandsOff} which ensures that the running cost is active in the first-order conditions and induces a nontrivial threshold in the Hamiltonian maximization. This creates a \emph{off} region in which actuation is not worthwhile. Abnormal extremals with \(\multip=0\) eliminate this threshold and may destroy sparsity. Assumption \ref{corr:ctrb} is a sufficient \emph{nondegeneracy} requirement for the ensemble \(\aset[]{A(\uparam)^k,b(\uparam)}_{k\ge 1}\): after averaging over \(\uparam\), the contributions from different parameter values cannot cancel in such a way that a nonzero terminal multiplier becomes invisible. In practical terms, any nontrivial multiplier inevitably produces a discernible change in the averaged switching signal, so the switching function cannot remain flat on a time interval. This rules out singular arcs and yields the bang–off–bang structure of the optimal control \(\cont\as(\cdot)\). The assumption is, in broad strokes, an ensemble-like counterpart of the usual controllability \cite{ref:UH:MS:EnsembleCond:2,ref:UH:MS:EnsembleCond:1, ref:dirrensemble} and without some form of it, averaging across parameters can mask the influence of the input --- even if each individual system is controllable --- so bang–off–bang cannot be guaranteed. The condition \ref{corr:ctrb} is verifiable and we refer the readers to \S\ref{sec:num_exp} for more details.
\end{remark}

\begin{proof}
Recall the switching function \( \lcrc{0}{\horizon} \ni t \mapsto \switch(t) \Let \int_{\pset} p(t; \uparam)^{\top} b(\uparam) \measure(\odif{\uparam}) \in \Rbb,\) and the maximization condition 
\begin{align}
\cont\as(t) \in  \argmax_{v \in \admcont} \aset[\big]{\switch(t)v - \eta \abs{v}},\nn
\end{align} from which we see that we see that
    \begin{align}
    \lcrc{0}{\horizon} \ni t \mapsto \cont\as(t) \in \begin{cases} \aset[]{-1} & \text{if } \switch(t) < - \multip \\ \lcrc{-1}{0} & \text{if } \switch(t) = - \multip \\  \aset[]{0} & \text{if } \abs{\switch(t)} < \multip\\ \lcrc{0}{1} & \text{if } \switch(t) = \multip \\ \aset[]{1} & \text{if } \switch(t) > \multip. \end{cases}    \nn 
\end{align}
Recall that \(\mu(\cdot)\) is the standard Eucledian Lebesgue measure; to show that \(\cont\as(\cdot)\) is bang-off-bang, i.e., \(\cont\as(t) \in \aset[]{-1,0,+1}\) for a.e. \(t \in \lcrc{0}{\horizon}\) we show that
\begin{align}
    \mu \bigl( \aset[]{t \in \lcrc{0}{\horizon} \suchthat \abs{\switch(t)} = \multip}\bigr) = 0,\nn
\end{align}
i.e., no singular arcs can exist.

Towards that end, we first establish that \(t \mapsto \switch(t)\) is an analytic function. Recall that, from \eqref{eq:adj:1:2}, the adjoint dynamics is given by \(-\dot{\adjx}(t;\uparam) = A(\uparam)^{\top} \adjx(\horizon;\uparam)\), and thus, for any fixed \(\uparam \in \pset\), we have \(\adjx(t;\uparam) = \epower{A(\uparam)^{\top}(\horizon - t)} p(\horizon;\uparam)\). Because \(\pset\) is compact and \(\uparam \mapsto (A(\uparam),b(\uparam))\) are continuous, one can find finite number \(M_A\) and \(M_b\) such that \(\norm{A(\uparam)} \le M_A\) and \(\norm{b(\uparam} \le M_b\) for all \(\uparam \in \actset\). Consequently, \(\sup_{\uparam,t}\norm{\epower{A(\uparam)^{\top}(T-t)}} < +\infty\), uniformly for all \((\uparam,t) \in \actset \times \lcrc{0}{\horizon}\). Using the power series expansion of the exponential, around some \(\hat{t} \in \lcrc{0}{\horizon}\), we have
\begin{align}
\switch(t)=\sum_{k=0}^{+\infty}\frac{(-1)^k}{k!}(t-\hat{t})^k  \int_{\actset} \inprod{A(\uparam)^k b(\uparam)}{\epower{A(\uparam)^{\top}(\horizon-\hat{t})}\adjx(\horizon;\uparam)}\actmeas(\odif{\uparam}). \nn 
\end{align} 
As the term inside the above integral is uniformly bounded and integrable, \(\adjx(\horizon;\cdot)\) is integrable (which follows from the transversality condition in the PMP), \(\actset\) is compact, and \(\actmeas\) is finite, the series converges absolutely and uniformly and thus 
\[
t \mapsto \switch(t)=\int_{\actset} \inprod{b(\uparam)}{\epower{A(\uparam)(\horizon-t)}\adjx(\horizon;\uparam) } \actmeas(\odif{\uparam}),
\]
is real-analytic. 

We are now ready to show the nonexistence of singular arcs. We proceed via contradiction: take a nonempty open interval \(J\subset \lcrc{0}{T}\) such that \(\abs{\switch(t)} =\multip\) for all \(t\in J\). Without loss of generality we consider the case \(\switch(t) = \multip\) on \(J\); the second case is identical with \(\multip \mapsto -\multip\). Real analyticity of \(\switch(\cdot)\) implies that it must be constant for all \(t \in \lcrc{0}{\horizon}\) \cite[Chapter 15]{ref:RudRealCom-1987}. Therefore for every \(k\ge1\), the \(k\)-th derivative \(\switch^{(k)}(t)=0\) for all \(t\in J\). To wit: 
\begin{align}\label{eq:switch:derivative:1}
\switch^{(k)}(t) & = (-1)^k\int_{\actset} \inprod{A(\uparam)^{k} b(\uparam)}{\epower{A(\uparam)^{\top}(\horizon-t)}\adjx(\horizon;\uparam)}
\actmeas(\odif{\uparam}).
\end{align}
Now, fix \(k\ge 1\) and choose any sequence \(t_n \uparrow \horizon\) as \(n\lra +\infty\). For \(\alpha \in  \actset\) we have 
\[
\inprod{A(\uparam)^{k} b(\uparam)}{\epower{A(\uparam)^{\top}(\horizon-t_n)}\adjx(\horizon;\uparam)} \to \inprod{A(\uparam))^{k} b(\uparam)}{\adjx(\horizon;\uparam)}
\]
as \(n \lra +\infty\), because \(\epower{A(\uparam)^{\top}(\horizon-t_n)} \to I\) and thus using the dominated convergence theorem
\begin{align}
    & \lim_{n \to +\infty} \int_{\actset} \inprod{A(\uparam)^{k} b(\uparam)}{\epower{A(\uparam)^{\top}(\horizon-t_n)}\adjx(\horizon;\uparam)} \actmeas(\odif{\uparam}) \nn \\& = \int_{\actset} \inprod{A(\uparam))^{k} b(\uparam)}{\adjx(\horizon;\uparam)} \teL \mathcal{M}_k.
\end{align}
Combining with \eqref{eq:switch:derivative:1} we get 
\begin{align}\label{eq:switch:derivative:2}
    \lim_{t \uparrow \horizon} \switch^k(t) = (-1)^k \mathcal{M}_k.
\end{align}
From Assumption \ref{corr:ctrb}, we have for any nonzero measurable \(v(\cdot)\) in the terminal normal cone, at least one \(M_{k}\) (with \(k \geq 1 \)) is nonzero. But from \eqref{eq:switch:derivative:2} we see that the singular interval enforces \(\mathcal{M}_1 = \mathcal{M}_2 = \ldots = \mathcal{M}_{d-1} =0\), hence \(v(\cdot) \Let \adjx(\horizon;\cdot) \equiv 0 \) \(\actmeas\)-a.e., and this gives us and \(\switch(t)\equiv 0 \) on \(\loro{0}{\horizon}\) --- a contradiction. Hence no singular interval \(J\) exists. Thus, \(\cont\as(\cdot)\) is bang-off-bang a.e., completing our proof.  
\end{proof}

We are ready to show the equivalence between the \(\lpL[0]\) problem \eqref{e:OCP} and the \(\lpL[1]\) problem \eqref{e:OCP:L1}. We recall and introduce some notations first. Recall that the set of admissible controls is given by 
\begin{align}
    \admu \Let \left\{\cont(\cdot) \in \lpL[\infty](\lcrc{0}{\horizon};\Rbb) \;\middle\vert\;  
    \begin{array}{@{}l@{}}
        \cont(t) \in \admcont \text{ for a.e. }t \in \lcrc{0}{\horizon} \text{ and }\\
        \psi\bigl(\st(\horizon;\uparam)\bigr) \le 0 \text{ for all }\uparam \in \pset
        \end{array}
        \right\}, \nn
\end{align}
which is, by assumption, nonempty. Also recall the expression of \(\lpL[0]\) objective function 
\begin{align}\label{eq:objlzero}
 \ \cont(\cdot) \mapsto \objlzero\bigl(\cont(\cdot)\bigr) \Let \mu \bigl( \aset[\big]{t \in \lcrc{0}{\horizon} \suchthat \abs{\cont(t)} \neq 0} \bigr),
\end{align}
and the \(\lpL[1]\) objective function
\begin{align}\label{eq:objlone}
 \cont(\cdot) \mapsto \objlone \bigl(\cont(\cdot)\bigr) \Let \int_{0}^{\horizon} \abs{\cont(t)} \odif{t}. 
\end{align}
We are ready to state our equivalence result. 

\begin{myOCP}
\begin{theorem}\label{thrm:equivalence}
Consider the OCP \eqref{e:OCP} and \eqref{e:OCP:L1} along with their associated data \ref{ocp:data:1}--\ref{ocp:data:3}. Recall the definition of the admissible set \eqref{eq:adm_cont_set} and objectives \eqref{eq:objlzero}--\eqref{eq:objlone}. Suppose that the hypotheses of Theorem \ref{thrm:L1:existence}, Theorem \ref{thrm:PMP:conditions}, and Corollary \ref{corr:bang-off-bang} hold. Define the optimal control set for the \(\lpL[0]\) problem by
\begin{align}\label{eq:optconlzero}
\optconlzero \Let \aset[\bigg]{ \cont_0\as(\cdot) \in \admu \suchthat \objlzero \bigl(\cont_0\as(\cdot) \bigr) \le \objlzero \bigl(\cont(\cdot) \bigr) \text{ for all }\cont(\cdot) \in \admu},
\end{align}
and the \(\lpL[1]\) problem by
\begin{align}\label{eq:optconlone}
\optconlone \Let \aset[\bigg]{ \cont_1\as(\cdot) \in \admu \suchthat \objlone \bigl(\cont_1\as(\cdot) \bigr) \le \objlone \bigl(\cont(\cdot) \bigr) \text{ for all }\cont(\cdot) \in \admu}.
\end{align}
Then the set of \(\lpL[0]\) optimal controls is equivalent to the set of \(\lpL[1]\) optimal controls, i.e., \(\optconlzero = \optconlone\).
\end{theorem}
\end{myOCP}


\begin{proof}
We first show \(\optconlone \subset \optconlzero\). Towards this end, let \(\cont_1\as(\cdot) \in \optconlone\); such \(\cont_1\as(\cdot)\) exists because due to Theorem \ref{thrm:L1:existence}. From Corollary \ref{corr:bang-off-bang} we know that \(t\mapsto \cont_1\as(t)\) must be a bang-off-bang, i.e., for a.e. \(t \in \lcrc{0}{\horizon}\), \(\cont_1\as(t)\in \aset[]{-1, 0, 1}\), which implies \(\abs{\cont_1\as(t)}\) is either \(0\) or \(1\). Thus, we have the equality
\begin{align}\label{eq:equality:1}
    \objlone \bigl( \cont_1\as(\cdot) \bigr) = \int_{0}^{\horizon}\abs{\cont_1\as(t)} \odif{t} = \int_{\supp(\cont_1\as)}1 \odif{t} = \mu \bigl( \supp(\cont_1\as) \bigr) = \objlzero \bigl( \cont_1\as(\cdot)\bigr) 
\end{align}
Now let \(\cont(\cdot) \in \admu\) be any feasible control trajectory. Then by design
\begin{align}\label{eq:inequality:2}
    \objlone \bigl(\cont(\cdot) \bigr) = \int_{0}^{\horizon}\abs{\cont(t)}\odif{t} = \int_{\supp(\cont)} \abs{\cont(t)}\odif{t} \le \int_{\supp(u)} 1 \odif{t} = \objlzero \bigl( \cont(\cdot) \bigr).
\end{align}
Thus, we have the chain
\begin{align}\label{eq:inequality:3}
    \objlzero \bigl( \cont_1\as(\cdot) \bigr) \stackrel{\mathclap{(\dag)}}{=} \objlone \bigl( \cont_1\as(\cdot) \bigr) \stackrel{\mathclap{(\ddag)}}{\le} \objlone \bigl(\cont(\cdot)\bigr) \stackrel{\mathclap{(\circ)}}{\le} \objlzero\bigl(\cont(\cdot)\bigr),
\end{align}
where the equality \((\dag)\) follows from \eqref{eq:equality:1}, the inequality \((\ddag)\) follows because \(\cont_1\as(\cdot)\) is \(\lpL[1]\)-optimal, and the final inequality \((\circ)\) follows from \eqref{eq:inequality:2}. Thus, \(\objlzero \bigl(\cont_1\as(\cdot) \bigr) \le \objlzero \bigl( \cont(\cdot) \bigr)\) for any feasible control \(\cont(\cdot)\), therefore \(\cont_1\as(\cdot) \in \optconlzero\).

We now prove the other direction \(\optconlzero \subset \optconlone\). Let \(\cont_0\as(\cdot) \in \optconlzero\) and also choose independently \(\cont_1\as(\cdot) \in \optconlone\). Because \(\cont_1\as(\cdot)\) is bang-off-bang, we have 
\begin{align}
    \objlone \bigl( \cont_1\as(\cdot) \bigr) = \int_0^{\horizon} \abs{\cont_1\as(t)} \odif{t} = \int_{\supp(\cont_1\as)} 1 \odif{t} = \objlzero \bigl( \cont_1\as(\cdot) \bigr),
\end{align}
and thus using \eqref{eq:equality:1} we have 
\begin{align}\label{eq:inequality:4}
    \objlzero\bigl( \cont_1\as(\cdot) \bigr) = \objlone \bigl(\cont_1\as(\cdot) \bigr) \le \objlone\bigl(\cont_0\as(\cdot) \bigr).
\end{align}
Now using the fact that \(\objlone \bigl(\cont_0\as(\cdot)\bigr) \le \objlzero\bigl(\cont_0\as(\cdot)\bigr)\) and that \(\cont_0\as(\cdot)\) is \(\lpL[0]\)-optimal, we have the inequalities
\begin{align}\label{eq:inequality:5}
\objlone \bigl(\cont_0\as(\cdot)\bigr) \le \objlzero\bigl(\cont_0\as(\cdot)\bigr) \le \objlzero\bigl( \cont_1\as(\cdot) \bigr).
\end{align}
Then, with \eqref{eq:inequality:4} and \eqref{eq:inequality:5} we have the chain
\begin{align}
      \objlzero\bigl( \cont_1\as(\cdot) \bigr) = \objlone \bigl(\cont_1\as(\cdot) \bigr) \le  \objlone \bigl(\cont_0\as(\cdot)\bigr) \le \objlzero\bigl(\cont_0\as(\cdot)\bigr) \le \objlzero\bigl( \cont_1\as(\cdot) \bigr),
\end{align}
which gives us \(\objlone\bigl(\cont_1\as(\cdot)\bigr) = \objlone\bigl(\cont_0\as(\cdot)\bigr)\) and thus \(\cont_0\as(\cdot) \in \optconlone\). The proof is complete. 
\end{proof} 

\begin{remark}
Theorem \ref{thrm:equivalence} states that, under normality and the nondegeneracy condition of Corollary \ref{corr:bang-off-bang}, the optimal solution sets of the robust \(\lpL[0]\) and robust \(\lpL[1]\) problems coincide, i.e., \(\optconlone = \optconlzero\). The practical implication is that the intrinsically nonconvex and nonsmooth robust \(\lpL[0]\) optimal control problems can instead be solved by solving their convex robust \(\lpL[1]\) surrogate, without loss of optimality, thereby enabling the use of a wide range of existing robust optimization solvers for numerical implementation; that's the subject matter of \S\ref{subsec:algo:dev:sip}.
\end{remark}

\subsection{Algorithmic developments}\label{subsec:algo:dev:sip}
Leveraging the equivalence result in Theorem \ref{thrm:equivalence}, we follow the scientific suggestions given in \cite[\S 1. Introduction, p.3]{ref:minmax:ness:vinter} to develop an algorithmic framework to solve the \(\lpL[1]\) robust OCP \eqref{e:OCP:L1}. OCP \eqref{e:OCP:L1} being an infinite-dimensional optimization problem over \(\admu\) is numerically intractable in general, and towards this end, we parametrize the space of control actions \(\admcon\) via piecewise constant dictionaries, which is a standard approach in numerical optimal control \cite{ref:DasGanAnjCha23CDC,ref:gansilcha24, ref:GanchaLCSS, ref:SG:SD:AA:MH:DC:sparse:num}.

Let \(N \in \N\) and consider a dictionary \(\dict \Let \aset[]{\reg_i(\cdot)}_{i \in \N} \subset  \pwcfunc_h(\lcrc{0}{\horizon};\Rbb)\) consisting of piecewise constant functions that are linearly independent and normalized, i.e., \(\max_{i} \max_{t}\abs{\reg_i(t)} =1.\) We define the parametrized set of admissible control trajectories \(\admconfunc_{\dict}  \subset \pwcfunc_h(\lcrc{0}{\horizon};\Rbb)\) by 
\begin{align*}
    \admconfunc_{\dict} \Let \linspan \aset[\big]{ \reg_i:\lcrc{0}{\horizon} \lra \Rbb \suchthat i = 1, \ldots, N }.
\end{align*}  
Let \(t \mapsto \Reg(t) \Let \bigl(\reg_1(t)\; \reg_2(t)\;\ldots \;\reg_N(t) \bigr) \in \Rbb^{N}\). Then, we define the parametrized control trajectory by
\begin{align}\label{eq:disct:control}
\lcrc{\tinit}{\tfin} \ni t \mapsto \cont^{\dict}(t) = \sum_{j=1}^{N} \Param_{j} \psi_j(t) = \Param\Reg(t)
\end{align}
where \(\Param \in \Rbb^{1 \times N}\) are the control parameters. Using the parametrizations in \eqref{eq:disct:control} we express the solution to \eqref{eq:sys} as
\begin{align}\label{eq:p_sol}
\lcrc{0}{\horizon} \ni t \mapsto \st_{\Param}(t; \param, \uparam) &\Let \epower{A(\uparam)t} \overline{x} + \int_{0}^{t} \epower{A(\uparam)(t-\tau)} b(\uparam) \Param \Reg(\tau)\odif{\tau}\text{ for }\uparam\in \pset,
\end{align}
and abusing notation we will write the trajectory \eqref{eq:p_sol} by \(t\mapsto \st_{\Param}(t;\uparam)\). The generating functions \(\reg_j(\cdot)\) for \(j=1,\ldots,N\) are piecewise constant and thus
\begin{align}    
  \Rbb^{N} \ni  \Reg(t) = \begin{cases}
        \Reg_{k-1} \quad &\text{if }t\in \lcro{t_{k-1}}{t_k},\\
        \Reg_K & \text{if }t\in \lcrc{t_K}{T},
    \end{cases}
\end{align}
where \(\Reg_{k-1}\) for \(k=1,\ldots,K-1\) and \(\Reg_{K}\) are all vectors in \(\Rbb^{N}\). The cost in \eqref{e:OCP:L1} can therefore, be further simplified as
\begin{align}\label{eq:simplified_cost}
    \objlone\bigl(\Param\Reg(\cdot)\bigr) =  &\Let \int_{0}^{\horizon} \abs{\Param \Reg(t)} \odif{t} 
    = \sum_{k=1}^{K} \int_{t_{k-1}}^{t_k} \abs{\Param \Reg_{k-1}}\odif{t} + \int_{t_K}^{\horizon}\abs{\Param\Reg_{K}}\odif{t}  \nn \\
    &  =\sum_{k=1}^{K}h\abs{\Param\Reg_{k-1}}+ (\horizon-t_K) \abs{\Param\Reg_K},
\end{align}
where \(h=t_k-t_{k-1}>0\).   
Putting everything together in \eqref{e:OCP:L1}, we obtain the following optimization problem:
 \begin{equation}
\label{eq:SR:pre_SIP} 
\begin{aligned}
& \hspace{-4mm}\inf_{\Param} 	&&   \objlone\bigl(\Param\Reg(\cdot)\bigr) = \sum_{k=1}^{K}h\abs{\Param\Reg_{k-1}} + (\horizon-t_K) \abs{\Param\Reg_K} \\
&  \hspace{-4mm}\sbjto		&&  \begin{cases}
\Param\Reg(t)\in \admcont, \psi\bigl(\st_{\Param}(\horizon;\uparam)\bigr) \le 0  \text{ for all }\uparam \in \pset.
\end{cases}
\end{aligned}
\end{equation}      
\begin{remark}
The optimal control problem \eqref{eq:SR:pre_SIP} has been reformulated as a finite-dimensional optimization problem. However, it still involves a \emph{compact, uncountable family of constraints} indexed by \(\uparam \in \pset\). Standard discretize-then-optimize methods for optimal control \cite{ref:betts-book,ref:SG:NR:DC:RB-22,ref:gansilcha24} are not suitable here, as they are unable to accommodate such an uncountable set of uncertainties. In the sequel, we develop and employ a robust optimization-based algorithm specifically designed to handle this class of problems.    
\end{remark}

The next proposition shows that \eqref{eq:SR:pre_SIP} has nice qualitative properties. 
\begin{myOCP}
\begin{proposition}\label{prop:unique solution}
Consider the OCP \eqref{e:OCP:L1} along with its problem data \ref{ocp:data:1}--\ref{ocp:data:3}, the parametrization \eqref{eq:disct:control}, and the OCP \eqref{eq:SR:pre_SIP}. Then, the set of admissible control parameters corresponding to the representation \eqref{eq:disct:control}
\begin{equation}\label{eq:ad_cont_traj}
\adparam \Let   \left\{\Param \in \Rbb^{1\times N} \, \middle\vert 
\begin{array}{@{}l@{}}
     \, \Param \Reg(t) \in \admcont \; \text{for all }t \in \lcrc{\tinit}{\tfin}
        \end{array}
        \right\},\nn
\end{equation}
is compact and convex, moreover, if the feasible set of \eqref{eq:SR:pre_SIP} is defined by
\begin{align}
\mathcal{F}_{O} \Let \left\{ \Param \in \adparam \;\middle\vert\;  
\begin{array}{@{}l@{}}
\psi \bigl(\st_{\Param}(\horizon;\uparam)\bigr)\le 0
\text{ for all }\uparam \in \pset
\end{array}
\right\},\nn
\end{align}
the program \eqref{eq:SR:pre_SIP} admits a solution.
\end{proposition}
\end{myOCP}

\begin{proof}
The mapping \(\Param \mapsto \Param\Reg(t)\) is linear and \(\admcont\) is convex; thus convexity of \(\adparam\) follows readily. Moreover, \(\admcont\) is closed and bounded and the mapping \(\Param \mapsto \Param\Reg(t)\) is continuous, thus \(\adparam\) is closed because it is the preimage of a closed set under a continuous map. The boundedness of \(\adparam\) follows immediately because \(\admcont\) is bounded and \(\reg_i(\cdot)\) are linearly independent. Thus \(\adparam\) is compact and convex. 

The set \(\mathcal{F}_{O}\) is compact. Indeed, \(\adparam\) is compact and for a fixed \(\uparam \in \pset\) (which is convex and compact), the mapping \(\adparam \ni \Param \mapsto \psi \circ \st_{\Param}(\horizon;\uparam)\) is continuous being composition of continuous maps \(\xi \mapsto \psi(\xi)\) and \(\adparam \ni \Param \mapsto \st_{\Param}(\horizon;\uparam)\). Thus, \(\mathcal{F}_{O}\) is the intersection of \(\adparam\) with the preimage of \(\lorc{-\infty}{0}\) under \(\adparam \ni \Param \mapsto \psi \circ \st_{\Param}(\horizon;\uparam)\) 
\begin{align}
    \mathcal{F}_{O} = \adparam \cap \bigcap_{\uparam \in \pset} \bigl(\psi \circ \st_{(\cdot)}(\horizon;\uparam)\bigr)^{-1}\bigl(\lorc{-\infty}{0}\bigr),
\end{align}
which is compact, being an intersection of a closed and a compact set \cite[Theorem 2.35, p.\ 37]{ref:Rud-Analysis}. The convexity of \(\mathcal{F}_O\) also follows immediately because \(\psi(\cdot)\) is convex and \(\adparam \ni \Param \mapsto \psi\bigl( \st_{\Theta}(\horizon;\uparam)\bigr)\) is affine. Finally, noting that \(\Param \mapsto \objlone\bigl(\Param\Reg(\cdot)\bigr)\) defined in \eqref{eq:simplified_cost} is continuous and convex, existence of solutions to \eqref{eq:SR:pre_SIP} follows readily from Weierstrass theorem \cite[Chapter 4, Theorem 4.16]{ref:Rud-Analysis}.
\end{proof}

We are ready to provide our main result concerning the SIP \eqref{eq:SR:pre_SIP} and its algorithmic solution whose proof follows from standard results in variational analysis \cite[Chapter 3]{ref:DonRoc-14}.
\begin{myOCP}
\begin{theorem}\label{thrm:exact:sol}
Consider the OCP \eqref{eq:SR:pre_SIP} and suppose that it is strictly feasible.\footnote{The OCP \eqref{eq:SR:pre_SIP}, satisfies a \emph{strict feasibility condition} if there exists a \(\Param\) such that \(\psi\bigl(\st(\horizon;\uparam)\bigr)<0\) for all \(\uparam \in \pset\).} Fix \(\param \in \Rbb^d\) and let the optimal value of \eqref{eq:SR:pre_SIP} be \(\valuefunc(\xz)\). Denote the set of feasible initial states \(\fsblset\) for \eqref{eq:SR:pre_SIP} by \(\fsblset \Let \aset[\big]{\param \in \Rbb^{d} \suchthat \valuefunc(\xz) < + \infty}\) and let \(\fsblset \neq \emptyset\). Let
\begin{align}\label{eq:dvar}
\dvar \Let  N = \text{dimension of the decision space of \eqref{eq:SR:pre_SIP}},
\end{align}
and define \(\pseq \Let (\uparam^1,\ldots,\uparam^{\dvar}) \in \pset^{\dvar}\). For a fixed \(\xz \in \fsblset\), define the relaxed version \(\pset^{\dvar} \ni \pseq \mapsto \gfunc(\xz;\pseq)\in \Rbb\), of \eqref{eq:SR:pre_SIP} by
 \begin{equation}
    \label{eq:g_func}
    \begin{aligned}
    \gfunc(\xz;\pseq) \Let   &\inf_{\Param\in \adparam}	&&  \objlone\bigl(\Param\Reg(\cdot)\bigr)  \\
        &  \sbjto && \begin{cases}
       \psi\bigl(\st_{\theta}(\horizon;\uparam^i)\bigr) \le 0
       \text{ for }i=1,\ldots,\dvar.
        \end{cases}
    \end{aligned}
\end{equation}   
\label{thrm:exact:sols:param-var}
 Consider the optimization problem
 \begin{equation}\label{e:global_max_prob_param-var}
\sup_{\pseq} \aset[\big]{\gfunc(\param;\pseq) \suchthat \pseq \in \pset^{\dvar}}.
\end{equation}
Then \(\pset^{\dvar} \ni \pseq \mapsto \gfunc\bigl(\xz;\pseq\bigr)\) is continuous and there exists \(\pseq^{\ast}(\xz)\) that solves \eqref{e:global_max_prob_param-var}. Moreover, we have exact solutions, i.e., \(\valuefunc(\xz)=\gfunc\bigl(\xz;\pseq^{\ast}(\xz)\bigr)\).
\end{theorem}
\end{myOCP}


\begin{proof}
For a fixed \(\param \in \fsblset\), to show the continuity of \(\gfunc(\param;\cdot)\) we will employ some standard results from variational analysis \cite[Chapter 3]{ref:DonRoc-14}. To this end, notice that the feasible set mapping of \eqref{eq:g_func} is nonempty and bounded, which follows from Proposition \ref{prop:unique solution}. As the mapping \(\Param \mapsto \objlone\bigl(\Param \Reg(\cdot)\bigr)\) is convex and the constraint map \((\Param,\uparam) \mapsto \psi \bigl(\st_{\Param}(\horizon;\uparam)\bigr)\) is jointly continuous, it follows that the feasible set mapping of \eqref{eq:g_func} is continuous \cite[Chapter 3, Example III.B]{ref:DonRoc-14} in the sense of Painlev\'e-Kuratowski \cite[Chapter 3, \S 3.2]{ref:DonRoc-14}. The preceding two properties imply that the mapping is Pompeiu-Hausdorff continuous \cite[Chapter 3, Theorem 3B.3]{ref:DonRoc-14} at \(\overline{\pseq} \in \pset^{\dvar}\).  Finally, employing \cite[Chapter 3, Theorem 3B.5]{ref:DonRoc-14} we conclude that \(\gfunc(\param;\cdot)\) is continuous around the point \(\overline{\pseq}\) which was arbitrary, and thus, \(\gfunc(\param;\cdot)\) is continuous. 

Continuity of \(\pset^{\dvar} \ni \pseq \mapsto \gfunc(\param;\pseq)\) along with the fact that \(\pset^{\dvar}\) is compact immediately gives us existence of an optimizer \(\pseq\as(\param)\).

As the optimization problem \eqref{eq:g_func} is a strictly feasible convex program with a convex and continuous cost, convex and jointly continuous (in the pair \((\Param,\uparam)\)) constraints, compact and convex decision space, and a compact uncertainty set, applying \cite[Theorem 1]{ref:DasAraCheCha-22} we readily obtain \(\valuefunc(\param) = \gfunc\bigl(\param,\pseq\as(\param)\bigr)\) for any \(\param \in \fsblset\). The proof is complete. 
\end{proof}

\begin{figure}[h]
\begin{center}
\begin{tikzpicture}[scale=0.7,
  node distance=1cm and 1.5cm,
  box/.style={
    draw, 
    rounded corners,
    minimum width=4cm,    
    minimum height=2cm,   
    align=center,
    drop shadow,
    fill=#1!30,
    draw opacity=0.1
  }
]
\node[box=Bittersweet!50!White, fill opacity=0.8, text opacity=1] (B) {(ii): \(\glbopt({\param})\)\\ for \eqref{e:global_max_prob_param-var}};
\node[box=Orchid!50!White, fill opacity=0.8, text opacity=1, above=of B] (A) {(i): Input data \(\param \in \fsblset\)};
\node[box=CornflowerBlue!50!White, fill opacity=0.8, text opacity=1, right=of B] (C) {(iii): \(\mathcal{P}\as(\param)\) \\ \text{parametric in }\param};
\node[box=Bittersweet!50!White, fill opacity=0.8, text opacity=1, below right=of B] (D) {(iv): OCP \eqref{eq:g_func};\\ parametric in \(\param\)};
\node[box=green!50!White, fill opacity=0.8, text opacity=1, below=of B] (E) {(v): \(t \mapsto \cont\as(t) = \Param\as(\param) \Reg(t)\)};
\node[draw=black, dashed, rounded corners, inner sep=10pt, fit=(B)] (Bbox) {};
\node[draw=black, dashed, rounded corners, inner sep=10pt, fit=(D)] (Bbox) {};
\draw[->,  line width=2pt] (A) -- (B);
\draw[->, line width=2pt] (B) -- (C);
\draw[->, line width=2pt] (C) -- (D);
\draw[->, line width=2pt] (D) -- (E);
\end{tikzpicture} 
\end{center}
\caption{A bird's-eye view of our algorithmic architecture \(\glbopt(\param)\).}
\label{fig:flowchat_ctmpc}
\end{figure}
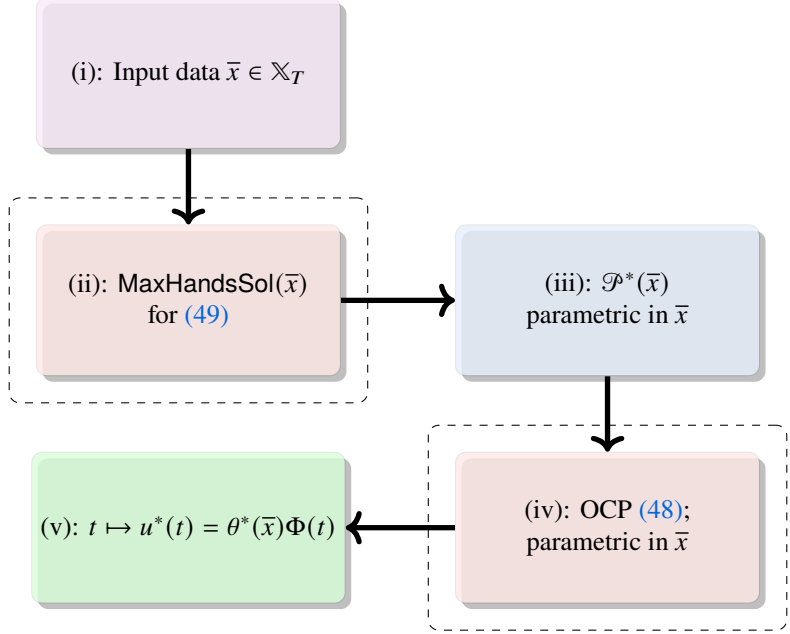

\begin{remark}
Theorem \ref{thrm:exact:sol} forms the fundamental backbone of our algorithmic approach to solving \eqref{e:OCP:L1}. It introduces a relaxed optimal control problem \eqref{eq:g_func} defined over a finite set of uncertainty realizations \(\pseq \Let (\uparam^1,\ldots,\uparam^{\dvar}) \in \pset^{\dvar}\), in contrast to the infinite family considered in \eqref{eq:SR:pre_SIP}. The first assertion, regarding the continuity of \(\gfunc(\param;\cdot)\), is crucial from a numerical standpoint, as it enables the use of numerous existing algorithms designed for global optimization problems with continuous objectives, such as \eqref{e:global_max_prob_param-var}. The second assertion guarantees the existence of at least one solution to problem \eqref{e:global_max_prob_param-var}. Most importantly, the final assertion precisely states that, by selecting \(\dvar\) uncertainty realizations and solving \eqref{e:global_max_prob_param-var} globally over \(\pset^{\dvar}\), the optimal value \(\gfunc\bigl(\param;\pseq\as(\param)\bigr)\) for a fixed \(\param \in \fsblset\) matches exactly the optimal value \(\valuefunc(\param)\) of the original semi-infinite program \eqref{eq:SR:pre_SIP}. A wide range of global optimization algorithms are readily available to solve \eqref{e:global_max_prob_param-var} efficiently, and our framework is algorithm-agnostic, provided that the chosen method complies with the conditions of Theorem \ref{thrm:exact:sol}; see the architecture \(\glbopt(\param)\) for details.\footnote{A few well-known algorithms are the Simulated annealing \cite{ref:CJPB-92}, differential evolution \cite{ref:storn1997differential}, SequOOL (Sequential Optimistic Optimization with Levels) \cite{ref:BarGabVal-19}, \cite{ref:GriValMun-15}, LIPO (Lipschitz Optimization based on Local Partitions) \cite{ref:MalVay-17}. Also see Fig. \ref{fig:flowchat_ctmpc} for a schematic.}
\end{remark}


{
\renewcommand{\algorithmcfname}{\(\glbopt(\param)\)}
\renewcommand{\thealgocf}{}
\begin{algorithm2e}[!ht]
\DontPrintSemicolon
\SetKwInOut{ini}{Initialize}
\SetKwInOut{giv}{Data}
\giv{Stopping criterion $\text{SC}(\cdot)$, threshold for the stopping criterion \(\tau\), fix \(\param \in \fsblset.\)}
\ini{initialize the constraint indices \( \bigl( \uparam^{\text{in},i}\bigr)_{i = 1}^{\dvar} \in \pset^{\dvar}\), initial guess for \(\gfunc_{\text{max}}\), initial guess for the solution \(\overline{\Param}\)}
   
\While{$\text{SC}(m) \leqslant \tau$}
{

Sample the uncertainty set \(\bigl(\uparam^{m,i}\bigr)_{i = 1}^{\dvar} \in \pset^{\dvar}\)
           
Solve the \emph{inner-problem} and evaluate \(\gfunc_m = \gfunc\Bigl(\bigl(\uparam^{m,i}\bigr)_{i = 1}^{\dvar};\param\Bigr)\) as defined in \eqref{eq:g_func}

\emph{Recover} the solution  \(\Param^m \in \argmin_{\widetilde{\Param}} \aset[\bigg]{\objlone\bigl(\widetilde{\Param}\Reg(\cdot)\bigr) \suchthat \text{ constraints in } \eqref{eq:g_func} \text{ hold at }  \bigl(\uparam^{m,i}\bigr)_{i = 1}^{\dvar} \in \pset^{\dvar}} \)

Update \( (\gfunc_{\text{max}}, \overline{\Param})\) \(\longleftarrow\) {IC}\((\gfunc_m, \Param^m)\) \,\,   (\({IC}\) is an improvement rule depending upon the global optimization routine)

Solve for \(\Param^m\) via any convex solver;
            
Update \(m \gets m+1\) \;
}

\caption{A general architecture to solve \eqref{e:global_max_prob_param-var}}
\label{alg:sip_algo}
\end{algorithm2e}
}

\begin{remark}
We note that classical methods to solve SIPs such as \cite{ref:blankenship1976infinitely} and methods in the spirit of Blankenship–Falk solve a sequence of finite relaxations by iteratively adding most violated constraints, so exactness is only achieved asymptotically (or depends on discretization/tolerances). In practice one must truncate after finitely many iterations, and these schemes typically provide no explicit bound on the error induced by truncation. In contrast, our framework yields an exact solution with finite --- and reasonably sized --- memory and computation: by Borwein’s result \cite[Theorem 4.1]{borwein1981direct} and \cite{ref:DasAraCheCha-22}, only \(\dvar\) samples of the uncertainty are required for exactness, and our algorithm returns them by solving a global optimization problem over \(\dvar\) uncertainty variables. As a result, no iterative constraint-generation loop or guess-and-refine sampling is needed, and no truncation/approximation error is incurred.
\end{remark}

\section{Discussion and numerical experiment}\label{sec:num_exp}
This section showcases a numerical simulation highlighting the effectiveness of our framework \(\glbopt(\param)\) applied to an uncertain continuous-time benchmark spring-mass-damper system. The global optimization component in \(\glbopt(\param)\) was implemented using the simulated annealing algorithm \cite{ref:CJPB-92}. All computations were performed in Julia version 1.1.14 on a laptop featuring an AMD Ryzen 5 5600H CPU and 8 GB of RAM.\footnote{S.G. thanks Ashwin Aravind for his assistance with the numerical simulations. }

Consider the spring-mass-damper system with parametric uncertainties
\begin{align}
\label{eq:num:smd:par}
\dot{\st}(t)  =  \begin{pmatrix}0& 1 \\ 
-2+\uparam_1& 0.6+\uparam_2
\end{pmatrix}  \st(t) + \begin{pmatrix}
    0 \\ 1+\uparam_3
\end{pmatrix} u(t) \quad\text{for all }t \in \lcrc{0}{\horizon}.
\end{align}
Let \(\uparam \Let (\uparam_1,\uparam_2,\uparam_3)^{\top}\). We have the following problem data 
\begin{align}
    \label{eq:num:prob:par:data}
    \hspace{-2mm}\begin{cases}
       \horizon=5,\,\st(0) \Let (-1, -1)^{\top},\,u(t) \in \lcrc{-1}{1},\, \uparam \in \pset \Let \lcrc{-0.1}{0.1}^3 \\ 
       \st(\horizon;\uparam) \in C \Let \left\{(\dummyx_1,\dummyx_2)^{\top}\in \Rbb^2 \;\middle\vert\;\begin{array}{@{}l@{}} \dummyx_1\leq 0.1 \text{ and}\\ \dummyx_2\leq 0.1\end{array}
\right\} \text{ for all }\uparam \in \pset.
    \end{cases}
\end{align}
Note that the condition \ref{corr:ctrb} can be easily verified; indeed: let \(\xi_1 \mapsto g_1(\xi_1) \Let \xi_1 - 0.1 \le 0\) and \(\xi_2 \mapsto g_2(\xi_2) = \xi_2 -0.1 \le 0\) and note that \(\nabla g_1(\xi_1) = e_1\) and \(\nabla g_2(\xi_2) = e_2\), where \(e_i\) for \(i=1,2\) are the standard Euclidean basis vectors. Then for the active parameters, \(N_C(\xi) = \aset[]{\mu_1 e_1 + \mu_2 e_2 \suchthat \mu_1\ge 0,\,\mu_2\ge 0}\) and thus \(\adjx(T;\uparam) = \mu_1(\uparam) e_1 + \mu_2(\uparam) e_2\) for \(\actmeas\)-a.e. \(\uparam \in \actset\). Since \(d=2\), we focus on the integral \(\int_{\actset}\inprod{A(\uparam)b(\uparam)}{\adjx(\horizon;\uparam)} \actmeas(\odif{\uparam})\) which is 
\begin{align}\label{eq:int:ineq}
\int_{\actset}\inprod{A(\uparam)b(\uparam)}{\adjx(\horizon;\uparam)} \actmeas(\odif{\uparam}) & = \int_{\actset} \hspace{-3mm} \bigl( \mu_1 \inprod{A(\uparam)b(\uparam)}{e_1} + \mu_2 \inprod{A(\uparam)b(\uparam)}{e_2}\bigr) \actmeas(\odif{\uparam}) \nn \\& \ge  \int_{\actset} \bigl( 0.9 \mu_1(\uparam) + 0.45 \mu_2(\alpha) \bigr) \actmeas(\odif{\uparam})  \nn \\& \ge 0.45 \int_{\actset}  (\mu_1(\uparam) + \mu_2(\uparam)) \actmeas(\odif{\uparam}).
\end{align}
If \(\adjx(\horizon;\uparam) \neq 0\), then the right-hand side of the inequality \eqref{eq:int:ineq} is positive and hence the condition \ref{corr:ctrb} holds.

\begin{figure}
\begin{subfigure}{0.49\textwidth}
\includegraphics[width=\linewidth]{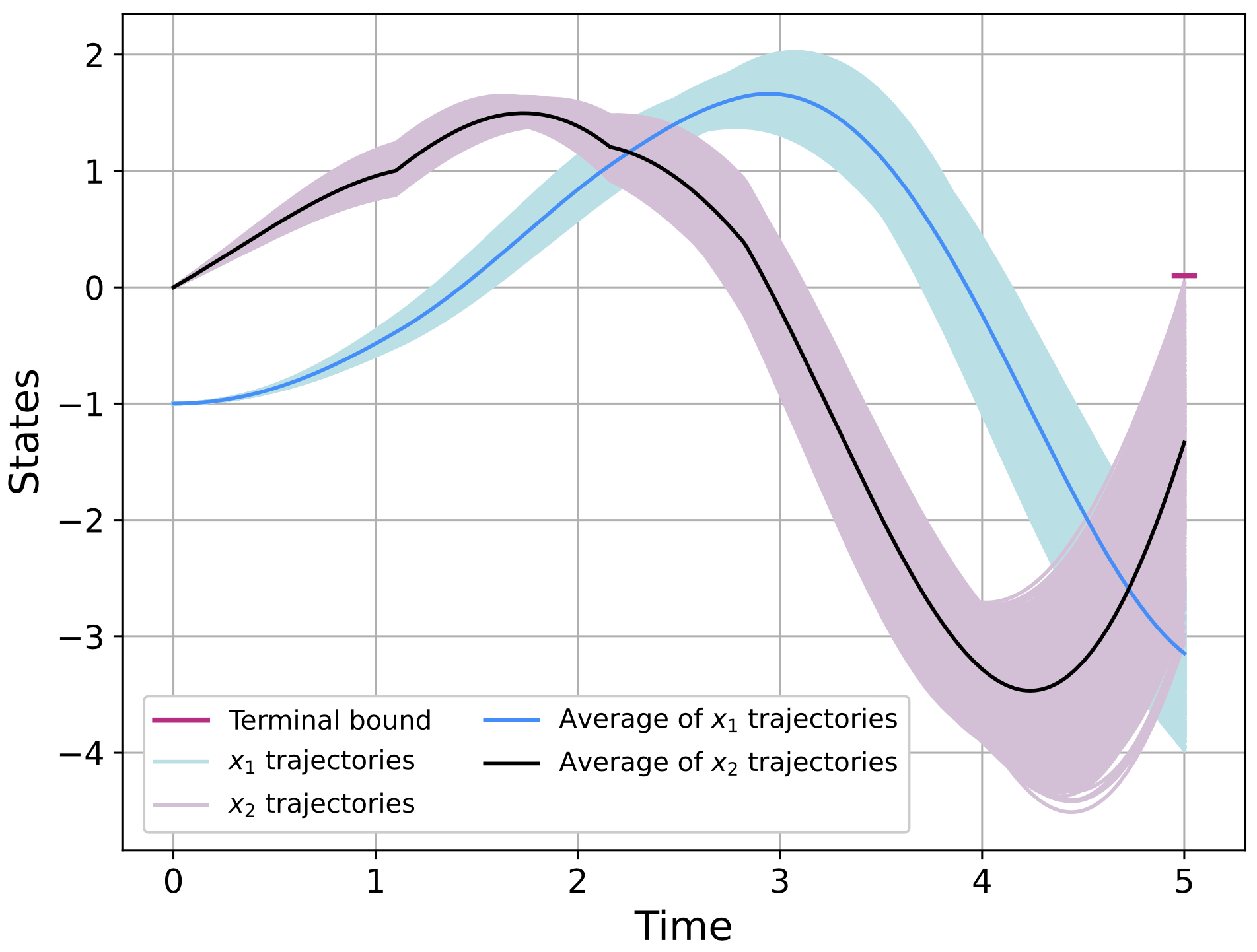}   
\caption{State trajectories.}
\label{fig:smd:par:traj}
\end{subfigure}
\begin{subfigure}{0.49\textwidth}
\includegraphics[width=\linewidth]{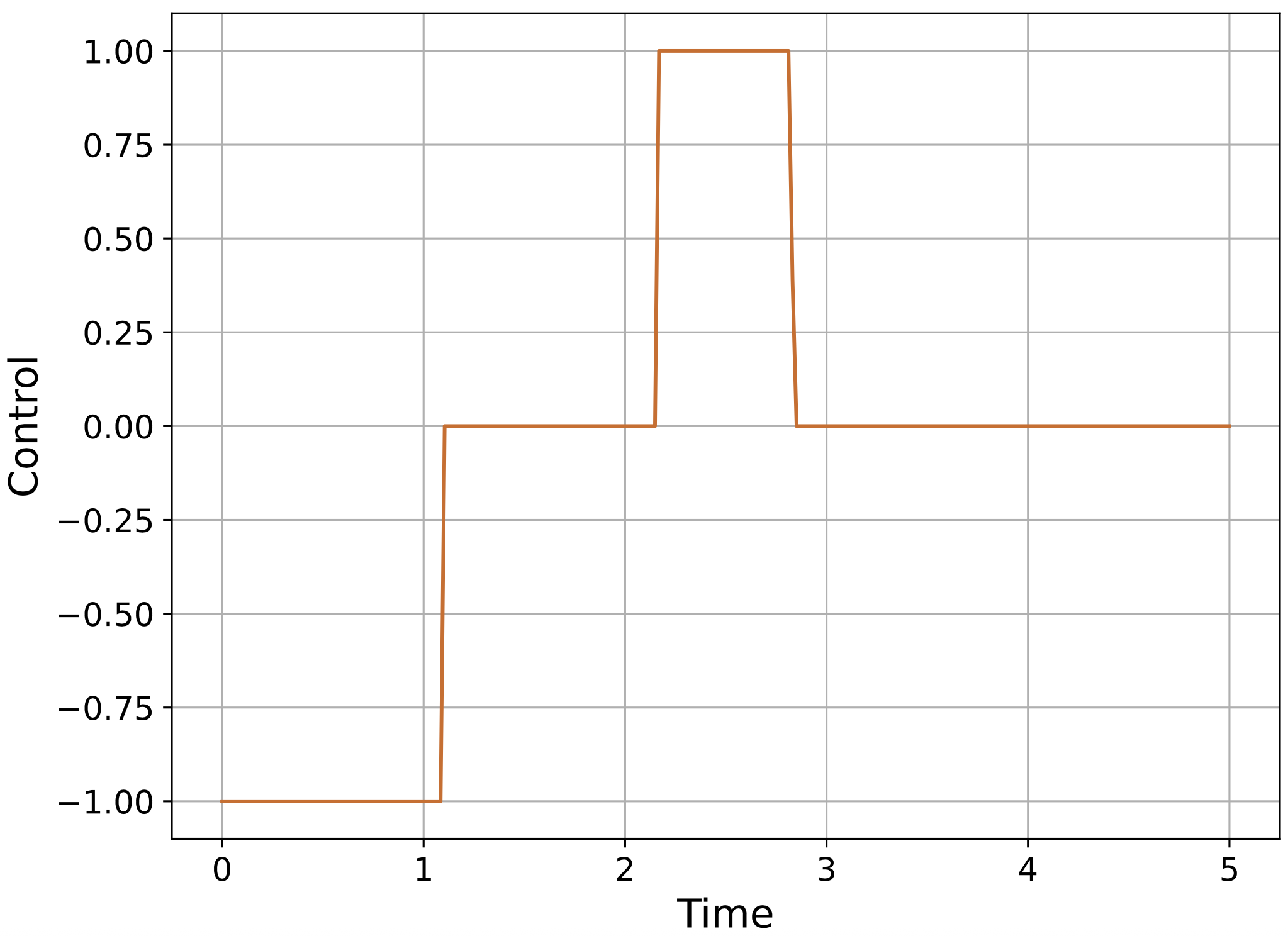}    
\caption{Control trajectory.}
\label{fig:smd:par:control}       
\end{subfigure}
\caption{State and control trajectories obtained via \(\glbopt\).}
\label{fig:smd:par}
\end{figure}

Now we focus on the numerical treatment. We parametrize the control trajectory \(\cont(\cdot)\) employing piecewise constant functions with \(N=250\). With this parametrization, we considered the OCP
\begin{align}
\label{eq:ocp:smd:par}
& \min_{\Param} && \objlone(\Param\Reg(\cdot)) \nn\\
& \sbjto && \begin{cases}
\text{dynamics }\eqref{eq:num:smd:par} \text{ and the data \eqref{eq:num:prob:par:data}}
\text{ for all }\uparam \in \pset,
\end{cases}
\end{align}
and employed the \(\glbopt(\cdot)\) architecture with the simulated annealing global optimization routine. The results of this numerical experiment are presented in Fig. \ref{fig:smd:par}. The state trajectories of the spring-mass-damper system \eqref{eq:num:smd:par}, subject to 10000 different realizations of the uncertain parameter, are shown in Fig.~\ref{fig:smd:par:traj}. The system is controlled using the sparse robust control computed by solving~\eqref{eq:ocp:smd:par} via \(\glbopt\bigl((-1, -1)^{\top}\bigr)\) architecture. The figure illustrates the effectiveness of the proposed approach in generating sparse robust control inputs for the system with data~\eqref{eq:num:prob:par:data}. As observed from Fig.~\ref{fig:smd:par:control}, the control input is sparse, i.e., being active only over a small portion of the time horizon. Despite the presence of parametric uncertainties, the resulting state trajectories satisfied both the state and terminal constraints, demonstrating the robustness of the proposed control scheme.

\begin{figure*}[t]
\centering
\begin{subfigure}{0.48\textwidth}
\includegraphics[scale=0.35]{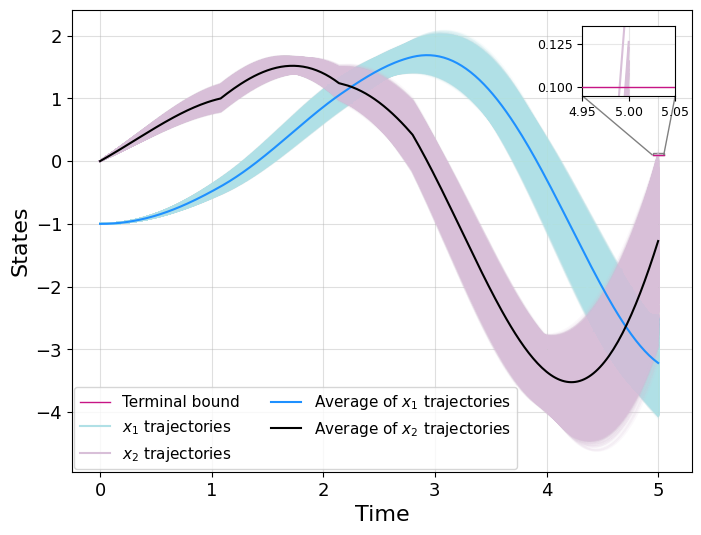}
\end{subfigure}\hfill
\begin{subfigure}{0.48\textwidth}
\includegraphics[scale=0.35]{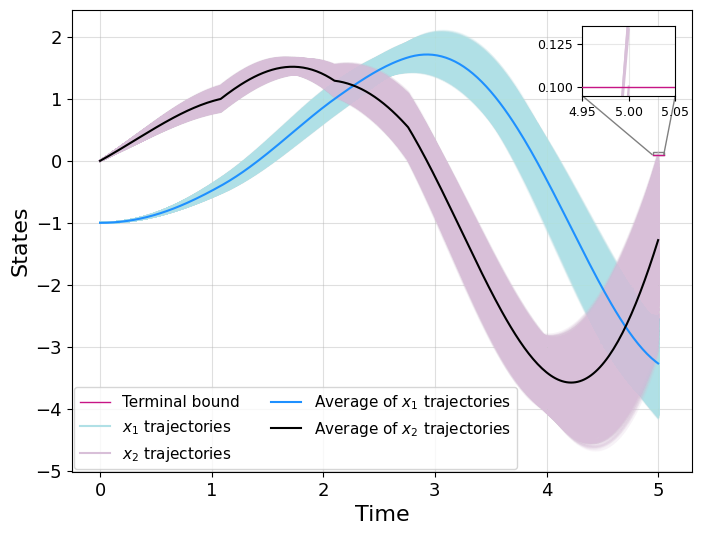}
\end{subfigure}
\caption{State trajectories obtained via the scenario approach, with \(N=500\) and \(1000\).}
\label{fig:smd:scenario}
\end{figure*}

For a baseline comparison, we present a few simulations based on the scenario robust optimization approach \cite{MC-SG:18}; an widely used technique in control and optimization. For the illustrations we considered 1000 and 5000 scenarios for the experiments. It can be clearly observed from Fig.~\ref{fig:smd:scenario} that multiple trajectories have violated the terminal constraints for all the scenarios presented here. Moreover, our algorithm, while employing finitary methods, achieves the true optimal values of CSIPs; see Table \ref{tab:metadata_ex_1}.
\begin{table}[t]
\centering
\begin{tabular}{lcc}
\toprule
Method  & value \(\big{(}\)for  \(\param = (-1,0)\)\(\big{)}\)
\\ \midrule
\(\glbopt(\param)\) & 89.01 \\
Scenario optimization \cite{MC-SG:18} (with \(N=200\)) & 82.93\\
Scenario optimization (with \(N=500\)) & 86.88  \\
Scenario optimization (with \(N=1000\)) & 87.01 \\
\bottomrule
\end{tabular}
\vspace{1mm}
\caption{Optimal values for \(\param \Let (-1\,\,0)^{\top}\), obtained via \(\glbopt(\param)\) and the scenario optimization algorithm.}
\label{tab:metadata_ex_1}
\end{table}

\section{Conclusion}
This article developed both a theoretical and an algorithmic framework for solving a class of robust sparse optimal control problems. On the theoretical front, we established an equivalence between the \(\lpL[0]\) and the \(\lpL[1]\) formulations of the robust optimal control problem. On the algorithmic side, we designed a robust optimization-driven architecture that solves such problems accurately and efficiently. Several promising avenues remain open for future investigation. One direction involves developing equivalence-type results and algorithms for more general robust problems, where both parametric uncertainties and external disturbances are present, including the design of faster gradient-based methods that leverage the regularity properties of \(\gfunc(\cdot)\).

\bibliographystyle{amsalpha}
\bibliography{refs}

\end{document}